\documentclass{article}
\usepackage{amsmath}
\usepackage{amsthm}
\usepackage{amssymb}
\usepackage{array}
\usepackage{titling}
\usepackage{hyperref}
\usepackage{cleveref}
\usepackage{mathtools}
\usepackage{biblatex}
\usepackage{titlesec}

\title{Exceptional Rays of Invariant Processes on Trees}
\author{Ori Gurel-Gurevich
\thanks{Einstein Institute of Mathematics, Hebrew University of Jerusalem. 
Email: \texttt{ori.gurel-gurevich@mail.huji.ac.il}.}
\and Nitzan Hen
\thanks{Einstein Institute of Mathematics, Hebrew University of Jerusalem. 
Email: \texttt{nitzan.hen@mail.huji.ac.il}.}
}
\predate{}
\postdate{}
\date{}

\titleformat*{\subsection}{\normalsize\bfseries}

\newtheorem{theorem}{Theorem}[section]
\newtheorem{proposition}[theorem]{Proposition}
\newtheorem{corollary}[theorem]{Corollary}
\newtheorem{claim}[theorem]{Claim}
\crefname{claim}{claim}{claims}
\Crefname{claim}{Claim}{Claims}

\newtheorem{lemma}[theorem]{Lemma}
\newtheorem{observation}[theorem]{Observation}
\crefname{observation}{observation}{observations}
\Crefname{observation}{Observation}{Observations}

\theoremstyle{definition}
\newtheorem{definition}[theorem]{Definition}
\newtheorem{fact}[theorem]{Fact}
\crefname{fact}{fact}{facts}
\Crefname{fact}{Fact}{Facts}

\theoremstyle{remark}
\newtheorem{remark}[theorem]{Remark}
\crefname{remark}{remark}{remarks}
\Crefname{remark}{Remark}{Remarks}

\newtheorem{example}[theorem]{Example}

\numberwithin{equation}{section}

\newcommand{\PP}{\mathbb{P}}
\newcommand{\EE}{\mathbb{E}}
\newcommand{\NN}{\mathbb{N}}
\newcommand{\ZZ}{\mathbb{Z}}

\newcommand{\RR}{\mathbb{R}}
\newcommand{\set}[1]{\left\{#1\right\}}
\newcommand{\one}{\mathbf{1}}
\newcommand{\eeps}{\varepsilon}

\newcommand{\MM}{\mathbf{M}}
\newcommand{\xx}{\mathbf{x}}

\newcommand{\lt}{<}
\newcommand{\gt}{>}

\DeclarePairedDelimiter{\floor}{\lfloor}{\rfloor}

\newcommand{\scloud}{s_{\text{cloud}}}
\newcommand{\sburst}{s_{\text{burst}}}
\newcommand{\ssust}{s_{\text{sust}}}

\DeclareMathOperator{\Aut}{Aut}
\DeclareMathOperator{\Vol}{Vol}
\DeclareMathOperator{\Var}{Var}
\DeclareMathOperator{\Cov}{Cov}

\DeclareMathOperator{\speed}{speed}

\newcommand{\gstar}{\mathcal{G}_\ast}
\newcommand{\gstarr}{\mathcal{G}_{\ast\ast}}
\newcommand{\pstar}{\mathcal{P}_\ast}

\newcommand{\as}{\text{a.s.}}

\begin{document}

\maketitle

\begin{abstract}
    Let \(A\) be an ergodic, non-constant random process on the \(d\)-regular tree \(T_d\), with marginal distribution \(X\). We prove that there always exist ``exceptional'' rays, whose asymptotic means are at least \(\EE(X)+\delta\); furthermore, for edge processes, \(\delta\gt 0\) can be taken to depend only on the law of \(X\). We compare different definitions of ``maximal asymptotic means'', relating them to different notions of speed of tree-indexed random walks. 
\end{abstract}

\tableofcontents

\clearpage

\section{Introduction}

Let \(T_d\) be the \(d\)-regular tree, and \(A\) be a random process on its edges or vertices. Suppose the distribution of \(A\) is invariant under all automorphisms of the tree, and let \(X\) denote the marginal distribution of a single edge or vertex. For every ray \(\gamma\) in \(T_d\), there exists an automorphism of \(T_d\) that acts as a shift on \(\gamma\); assuming that \(A\) is ergodic, this implies that the partial averages along \(\gamma\) converge a.s. to \(\EE(X)\) (for i.i.d. processes, this follows since \(A\) restricted to \(\gamma\) is again an i.i.d. process; for general ergodic processes this claim is not immediate, see \Cref{corollary:ray-ergodicity}).

However, since the set of such rays is uncountable, this does not rule out the existence (with positive probability) of ``exceptional'' rays, with asymptotic means larger or smaller than \(\EE(X)\). And indeed, this phenomenon has been thoroughly demonstrated for i.i.d. processes, and in particular for Bernoulli percolation.

In \cite{benjamini2018invariant}, Benjamini and Gurel-Gurevich show similar behavior for the class of automorphism-invariant percolation processes on the edges of \(T_d\). More precisely, given such a percolation process \(A\) they define \[
D(\gamma)=\limsup_{n\to\infty}\frac{1}{n}\sum_{k=1}^n A(\gamma_{k-1},\gamma_k),\quad D(A)=\sup_{\gamma} D(\gamma)
\] 
and \(D_d(p)=\inf_{A}D(A)\), where the infimum goes over all invariant bond percolation processes \(A\) on \(T_d\), with marginal probability \(p\) of an edge being open. They showed that \(D_d(p)\gt p\) for all \(0\lt p\lt 1\) and \(d\ge 4\), as well as for some values of \(p\) for \(d=3\). In other words, there are rays with densities not only larger than \(p\), but larger than \(p+\delta\), where \(\delta\) depends only on \(p\) and \(d\) (but not on the dependence between different edges, i.e. on their joint distribution).

Aside from the remaining values of \(p\) in the case \(d=3\), they also asked a few other natural questions: \begin{enumerate}
    \item Does the function \(D_d(p)\) remain the same if we replace the \(\limsup\) in \(D(\gamma)\) by a \(\liminf\)? If not, do the results still hold in that case?
    \item Considering more general invariant processes, are there always rays with asymptotic means larger than the expected value?
    \item Can a similar result be shown in the case of vertex percolation?
\end{enumerate}

We answer these questions below. It turns out that a natural framework for these questions of "maximal asymptotic means" is that of speeds of \textit{tree-indexed random walks}: given a tree \(T\) with root \(o\), put labels \((X_w)_{w\in V(T)}\), typically i.i.d., on the vertices. Then define \[
    S_{o,v}=\sum_{o\le w\le v}X_w
\]
where the sum goes over the unique path \([o,v]\) between \(o\) and \(v\). The process \((S_v)_{v\in V(T)}\) is called a \textbf{tree-indexed random walk}. When \(T=\NN\) (rooted at \(0\)) and \((X_i)_i\) are i.i.d. and uniform on \(\set{-1,1}\), we get the familiar simple random walk on \(\ZZ\); for a general tree, we may analogously view the walk as a collection of particles moving in space at increments \((X_w)_w\), however now a particle could also "split" into multiple descendants. 

In particular, since the walk now consists of many particles, the standard definition of speed – as the limit of \(\frac{d(S_n,S_0)}{n}\) – no longer makes sense. It is replaced by three other natural notions:

\begin{definition}[Speeds] \label{def:speeds}
Let \(A\) be a process on the vertices or edges of a tree \(T\) rooted at \(o\). Define its
\[
\setlength\arraycolsep{1pt}
\begin{array}{@{\hspace{1em}} p{4 cm}l>{\displaystyle}l p{2.5 cm}}
    \text{\textbf{cloud speed} by }&\scloud(A)&:= \limsup_{n\to \infty} \frac{1}{n}\max_{|v|=n}A_{o,v}, & \\[12pt]
    \text{\textbf{burst speed} by }&\sburst(A)&:= \sup_{\zeta\in \partial T}\limsup_{n\to\infty} \frac{1}{n}A_{\zeta_0,\zeta_n}, & \\[12pt]
    \text{\textbf{sustainable speed} by }\quad&\ssust(A)&:=\sup_{\zeta\in \partial T}\liminf_{n\to\infty} \frac{1}{n}A_{\zeta_0,\zeta_n}& \\
\end{array}
\]

Where \(A_{x,y}\) denotes the sum, with respect to \(A\), along the path \([x,y]\).
\end{definition}

Following \cite{benjamini1994tree}, we interpret the names as follows: for the cloud speed we treat the level set \(\set{v: |v|=n}\) as the "cloud" of particles at time \(n\), and measure the greatest displacement among them. In contrast, the sustainable and burst speed track individual lines of descent; they capture the greatest rate of escape along a single ray, either consistently (\(\ssust\)) or along some subsequence of levels (\(\sburst\)). \\
We always have \(\scloud\ge \sburst\ge \ssust\). In the i.i.d. case on a general tree, or the ergodic case on \(T_d\), these speeds are a.s. constant and \(\ssust(A)\ge \EE(X)\). Most results in the literature are concerned with i.i.d. processes, on a general tree \(T\) (see \cite[Sections 13.8,13.9]{lyonsperes}).

In the general setting of automorphism-invariant processes on the tree, the tree-indexed random walk interpretation of these notions is not as clear as in the i.i.d. case. For example, invariance implies that the process "looks the same" from any vertex, while the random walk interpretation clearly attributes special meaning to the root \(o\) and the levels of the tree determined by it. However, the definitions above still capture the essence of "maximal means" achieved by the process. Note, in particular, that the quantity \(D(A)\) from Banjamini and Gurel-Gurevich is exactly \(\sburst(A)\).

We now state our main results. First, since we work only on regular trees, one might ask whether these different speeds actually coincide (this is also referred to by question 1 above).

\begin{theorem} \label{thm:scloud-equals-sburst}
    Let \(A\) be an invariant process on the edges or vertices of \(T_d\). Then a.s. \[
        \scloud(A)=\sburst(A)
    \]
\end{theorem}

\begin{theorem} \label{thm:process-where-scloud-gt-ssust-exists}
    There exists an invariant process on the edges of \(T_d\) for which a.s. \[
        \scloud(A)\gt \ssust(A)
    \]
\end{theorem}

Next, we consider the existence of exceptional rays, with asymptotic means larger than \(\EE(X)\). We provide a positive answer to question 2 above, under the assumption that \(\Var(X)\lt\infty\), and prove even more: for ergodic processes, we can find a gap that depends only on the marginal distribution \(X\).

\begin{theorem} \label{thm:exceptional-rays-exist}
    Let \(A\) be an ergodic process on the edges of \(T_d\). Suppose the marginal distribution \(X\) of \(A\) is non-constant and has finite variance. Then a.s. \[
    \ssust(A)\ge \EE(X)+\delta_X
    \]
    where \(\delta_{X}\gt 0\) depends only on the law of \(X\).
\end{theorem}

Using ergodic decomposition, an immediate corollary of \Cref{thm:exceptional-rays-exist} is that any invariant (possibly non-ergodic) process \(A\) that satisfies \(\Var(X)\lt\infty\) and \(\EE(\ssust(A))=\EE(X)\) is a mixture of constant processes.

A central tool in the proof of \Cref{thm:exceptional-rays-exist} is a family of \textit{random walks on unimodular random networks} (RWRN), supported on trees (see \Cref{section:random-walk-unimodular-random-networks} for precise definitions). One of our results is determining the speed of such walks:

\begin{theorem}
    Let \(\xx\) be an ergodic, transient RWRN supported on a tree, and \(\xi=\xi(\xx)\)  its limiting end. Then a.s. \[
        \speed(\xx)=\EE_\nu(d(x_0,x_1))-2\EE_\nu(|x_{-1}\wedge \xi|)
    \]
    where \(|x_{-1}\wedge \xi|\) denotes the number of edges shared between the path \([x_0, x_{-1}]\) and the walk's escape ray \([x_0,\xi)\), and \(\EE_\nu\) denotes expectation with respect to the stationary measure of the walk.
\end{theorem}

Finally, we turn to vertex processes. As observed by Benjamini and Gurel-Gurevich, in this class there actually exists a natural, non-constant invariant percolation process for which a.s. \textit{all} rays have density \(p=\frac{1}{2}\): note that \(T_d\) is bipartite; let \((P,Q)\) be a partition. Define a percolation \(S\) such that \(S|_P\equiv 1,S|_{Q}=0\) or \(S|_{P}\equiv 0, S|_{Q}\equiv 1\), each with probability \(\frac{1}{2}\). Then \(S\) is ergodic, but a.s. \(\scloud(S)=\frac{1}{2}\) (and hence \(\ssust(S),\sburst(S)\) also equal \(\frac{1}{2}\) a.s.). \\
Of course, the same applies to the collection of processes \(\mathcal{S}_{a,b}\) defined the same way, but with values \(a,b\in \RR\) replacing \(0,1\). We show that these processes are effectively the only ones:

\begin{theorem} \label{thm:vertex-exceptional-rays-exist}
    Let \(A\) be an ergodic process on the vertices of \(T_d\), whose marginal distribution \(X\) has finite variance. If \(\ssust(A)=\EE(X)\), then \[
    A\overset{d}{=}\mathcal{S}_{a,b}
    \]
    for some \(a,b\in \RR\).
\end{theorem}

As an immediate result (via ergodic decomposition), any invariant process on the vertices with \(\EE(\ssust(A))=\EE(X)\) (and \(\Var(X)\lt\infty\)) is a mixture of the processes \(\set{\mathcal{S}_{a,b}}_{a,b}\). \\

The rest of the text is structured as follows. \Cref{section:preliminaries} is devoted to preliminaries, i.e. clarifying notation and providing basic definitions and results that are used throughout the paper. In \Cref{section:scloud-equals-sburst} we prove \Cref{thm:scloud-equals-sburst}, and in \Cref{section:worm-process} we construct a percolation process that proves \Cref{thm:process-where-scloud-gt-ssust-exists}. \Cref{section:random-walk-unimodular-random-networks} presents some general results about unimodular random networks and random walks on them, which we use in \Cref{section:existence-of-exceptional-rays} to prove Theorems \ref{thm:exceptional-rays-exist} and \ref{thm:vertex-exceptional-rays-exist}.

\section{Preliminaries} \label{section:preliminaries}

\subsection{Graphs and trees}

Let \(T=(V,E)\) be a tree, i.e. a connected graph without cycles. Denote by \(\deg_T x\) the degree of a vertex \(x\in V\). A \textbf{walk} on \(T\) is simply a finite or (countably) infinite sequence of vertices \(\mathbf{w}=(w_0,w_1,...)\). Note that we do not require that \(w_i\) is adjacent to \(w_{i-1}\) in the above; when successive vertices are indeed adjacent, we say that the walk is a \textbf{nearest-neighbor} walk, or that its \textbf{step size} is 1. A \textbf{path} is a self-avoiding nearest-neighbor walk, i.e. it does not pass through any vertex twice. We also call an infinite path a \textbf{ray}.
For vertices \(x,y\), denote by \([x,y]\) the unique path in \(T\) from \(x\) to \(y\); if \(z\) is another vertex, write \(x\le z\le y\) if \(z\) is on that path. \\
We will only be dealing with trees that are infinite and \textbf{locally finite}, which means that \(\deg_T x<\infty\) for every vertex \(x\). This implies, in particular, that the vertex and edge sets of the tree are countable. Our trees will typically be \textbf{rooted}, which simply means there's a distinguished vertex \(o\in V\) called the \textit{root}.

A central example, which will be of primary interest to us, is that of regular trees: for \(d\ge 3\), let \(T_d\) be the \textbf{\(d\)-regular tree}, that is the unique (infinite) tree for which every vertex has degree \(d\). \(T_d\) can be represented as the Cayley graph of the group \[
    G_d=\langle a_1,...,a_d\mid a_i^2=e\rangle
\]
with repsect to the generators \(a_1,...,a_d\). \\
Denote by \(d_T\) the graph metric on \(T\); since \(T\) is a tree, \(d_T(x,y)\) is simply the length of the path \([x,y]\). For a rooted tree, we also write \(|x|:=d_T(x,o)\). Note that \(d_T(x,z)=d_T(x,y)+d_T(y,z)\) if and only if \(x\le y\le z\). For a vertex \(x\) and a radius \(r\ge 0\), let \[
S_r(x)=\set{y\in V\mid d_T(x,y)=r},\quad B_r(x)=\set{y\in V\mid d_T(x,y)\le r}
\] be the \textbf{sphere} and \textbf{(closed) ball} of radius \(r\) around \(x\), respectively. This notation extends to sets: given subsets \(K, K'\subseteq T\), define \(d_T(x,K)=\min_{y\in K}d_T(x,y)\), \(d_T(K',K)=\min_{x\in K'}d_T(x,K)\) and \[
S_r(K)=\set{y\in V\mid d_T(y,K)=r},\quad B_r(K)=\set{y\in V\mid d_T(y,K)\le r}
\]
We call \(B_r(K)\) the \textbf{\(r\)-buffer} of \(K\). \\

An \textbf{isomorphism} of two graphs \(G, H\) is a bijection \(\varphi: V(G)\to V(H)\) such that \(\varphi(x)\sim \varphi(y)\) -- that is, \(\varphi(x), \varphi(y)\) are connected by an edge -- if and only if \(x\sim y\). An \textbf{automorphism} of \(G\) is an isomorphism from \(G\) to itself. The set of automorphisms of \(G\) is denoted \(\Aut(G)\); it is a group under composition, and naturally acts on \(G\). We will only be interested in automorphisms on \(T_d\), whose automorphism groups are very rich: for one, the action \(\Aut(T_d)\curvearrowright T_d\) is \textbf{transitive}, meaning that for every two vertices \(x,y\), there is an automorphism \(\varphi\in \Aut(T_d)\) that maps \(x\) onto \(y\). However a much stronger property holds:

\begin{fact} \label{fact:regular-tree-transitivity}
    Let \(G,H\) be two connected subgraphs of \(T_d\). If \(G, H\) are isomorphic as graphs, then there exists an automorphism \(\varphi\) of \(T_d\) that maps \(G\) onto \(H\).
\end{fact}

In fact, this property characterizes regular trees: suppose \(G\) is an infinite, locally finite, connected graph that contains a cycle of length \(\ell\). Then it contains a path \((v_0,v_1,...,v_{\ell-1})\) with \(d_G(v_0,v_{\ell-1})=1\). Since \(G\) is connected and infinite, it must also contain a ``geodesic'' path \((u_0,u_1,...,u_{\ell-1})\) with \(d_G(u_0,u_{\ell-1})=\ell-1\). These are both paths of length \(\ell-1\), but there cannot be an automorphism of \(G\) mapping one onto the other. Thus any graph satisfying the property at \Cref{fact:regular-tree-transitivity} must be a tree; such a graph must also be transitive, hence it must be a regular tree. \\

Another concept that will prove useful to us is that of an \textbf{end} of a tree. We say two rays \(\gamma_1,\gamma_2\) in a tree \(T\) are \textit{end-equivalent} if the symmetric difference \(V(\gamma_1)\triangle V(\gamma_2)\) is finite, that is if the rays share all but finitely many vertices. An \textbf{end} in \(T\) is an equivalence class of rays under this equivalence relation. We call the set of ends the \textbf{boundary} of \(T\), denoted \(\partial T\). Note that two end-equivalent rays with the same starting point necessarily coincide. In other words, a ray \(\gamma\) in \(T\) is determined by its starting point \(x=\gamma_0\) and its end \(\zeta=[\gamma]\) -- we then denote \(\gamma=[x,\zeta)\).

Ends can also be described from another perspective. Given three vertices \(x,y,z\), define \(y\wedge_x z\) to be the vertex in the intersection \([x,y]\cap [x,z]\) that's farthest from \(x\); since the intersection is a path, this means exactly that \[
    [x,y\wedge_x z]=[x,y]\cap [x,z]
\] 
Alternatively, \(y\wedge_x z\) is the unique vertex that belongs to all three paths \([x,y]\),\newline
\([x,z]\) and \([y,z]\) (and in particular, \(y\wedge_x z=z\wedge_y x=x\wedge_z y\)). A key identity relating the distances of these paths is \begin{equation} \label{eqn:wedge-distance-identity}
    d(x,y\wedge_x z)=\frac{1}{2}(d(x,y)+d(x,z)-d(y,z))   
\end{equation}
Now given a root \(o\), we write \(x\wedge y:=x\wedge _o y\). Define a metric \(d'\) on \(T\) by \begin{equation} \label{eqn:end-metric-on-T_d}
    d'(x,y):=\begin{cases}
        0 & x=y \\
        e^{-|x\wedge y|} & x\neq y
    \end{cases}   
\end{equation}
That is, vertices \(x,y\) are close under \(d'\) if they share a long initial segment. We then consider the completion \(\bar{T}\) of \(T\) with respect to this metric. The boundary \(\partial T\) may be redefined as \[
    \partial T:=\bar{T}\setminus T
\]
and an end as a point \(\zeta\in\partial T\). Observe that while the metric \(d'\) itself depends on the choice of root \(o\), the underlying topology, and in particular the boundary \(\partial T\), do not. \\
The operation \(x\wedge y\) naturally extends to \(\bar{T}\): given a vertex \(x\) and an end \(\zeta\), define \(x\wedge \zeta\) to be the vertex in the intersection \([o,x]\cap [o,\zeta)\) that's farthest from the root. \(\zeta_1\wedge \zeta_2\), for two ends \(\zeta_1,\zeta_2\), is defined similarly. \\
A walk \(\mathbf{w}\) on \(T\) is said to be \textbf{end-convergent} if there exists an end \(\zeta\) such that \(w_n\to \zeta\) under the metric \(d'\), that is if \(|w_n\wedge \zeta|\to \infty\). We then call \(\zeta\) the \textbf{limiting end} of the walk \(\mathbf{w}\), and denote it by \(\xi(\mathbf{w})\). We also call the canonical representative \(\gamma=[o,\xi(\mathbf{w}))\) the \textbf{escape ray} of the walk. Since \(\bar{T}\) is complete, \(\mathbf{w}\) is end-convergent if and only if it is a Cauchy sequence and is not eventually constant. End-convergence can also be given an equivalent geometric formulation: \(w_n\to \zeta\) if and only if for every representative ray \(\gamma=[x,\zeta)\) (equivalently for some such \(\gamma\)), and for every \(k\), there exists some \(N\) such that for all \(n\gt N\), \begin{equation} \label{eqn:end-convergence-geometric-condition}
    x\le \gamma_k\le w_n 
\end{equation}
A walk \(\mathbf{w}\) is called \textbf{transient} if it eventually escapes every finite set, i.e. if for every finite subgraph \(K\subseteq T\), \(w_n\notin K\) for all large enough \(n\). Every end-convergent walk is transient; the converse is not true in general, but a natural sufficient condition is that the walk has bounded step size, i.e. there exists some \(M\) such that \(d_T(w_i,w_{i+1})\le M\) for all \(i\).

\subsection{Random processes}
An  \textbf{\(I\)-indexed (random) process} is a collection of random variables \(A=\set{A_i}_{i\in I}\) defined on the same probability space. Alternatively, we may view \(A\) as a random element in \(\RR^{I}\). Given a tree \(T\), a random \textbf{edge} or \textbf{vertex process} is simply a random process on \(E(T)\) or \(V(T)\), respectively. In cases where the context is clear or the argument for edge and vertex processes is the same, we will use "random process" or just "process". We will mostly be interested in processes whose \textit{marginal distributions} -- that is, the distributions of each random variable \(A_i\) on its own -- are all identical. By a slight abuse of notation, we will usually refer to the marginal distribution as another, "abstract" random variable \(X\) whose distribution is the same as \(A_i\).
An important particular class of processes is percolation: a \textbf{(general) percolation} is simply a random process with values in \(\set{0,1}\). More specifically, a \textbf{bond percolation} is a \(\set{0,1}\)-valued edge process, and a \textbf{site percolation} is a \(\set{0,1}\)-valued vertex process. In this context we may identify the process with a random subgraph of \(T\), which consists of the edges or vertices with value \(1\); we call such an edge or vertex \textbf{open} (otherwise, i.e. when its value is \(0\), we call it \textbf{closed}). \\
Another class of processes, perhaps the canonical one, is that of \textbf{i.i.d.\ processes}, in which the random variables \(\set{A_i}_{i\in I}\) are independent and identically distributed. A particular case of interest is \textbf{Bernoulli percolation}, wherein each edge or vertex is open independently with some fixed probability \(p\) (or in other words, the marginal distribution is \(\operatorname{Bernoulli}(p)\)). \\
A random process \(A\) is called \textbf{\(m\)-dependent} if for every \(K_1,K_2,...,K_\ell\subseteq T\) with \(d_T(K_i,K_j)\gt m\) for all \(i\neq j\), the restrictions \(A|_{K_1},...,A|_{K_\ell}\) are independent. \(A\) is \textbf{finitely dependent} if it's \(m\)-dependent for some \(m\). \\

The class of processes that we will be most interested in, however, is invariant processes. Let \((\Omega,\mathcal{F})\) be a measurable space and \(\Gamma\) a group acting on it. A measure \(\mu\) on \(\Omega\) is called \textbf{\(\Gamma\)-invariant} if for every event \(E\) and every \(\gamma\in \Gamma\), \[
    \mu(\gamma.E)=\mu(E)
\]
An edge or vertex process \(A\) is called \textbf{(automorphism) invariant}, or \textbf{stationary}, if its law is an \(\Aut(T_d)\)-invariant measure (on \(\RR^{E(T_d)}\) or \(\RR^{V(T_d)}\), respectively).
An event \(E\) is called \textbf{invariant} if \(\phi.E=E\) for all \(\phi\in \Aut(T_d)\). An invariant measure \(\mu\) on \(\RR^{E(T_d)}\) or \(\RR^{V(T_d)}\) is called \textbf{ergodic} if for every invariant event \(E\), \(\mu(E)\in \set{0,1}\). 
Equivalently, \(\mu\) is ergodic if and only if every invariant function is a.e. constant, where we call a measurable function \(f: \Omega\to \RR\) invariant if for all \(\gamma\in\Gamma\), \(f(\gamma.\omega)=f(\omega)\) for almost every \(\omega\).
A process \(A\) is called ergodic if its law is an ergodic measure. A central result in ergodic theory is that every \(\Gamma\)-invariant measure on \(\Omega\) can be uniquely represented as a mixture of ergodic measures on the same space (see for example \cite[Theorem 3.22]{Glasner2003}). In particular, if \(F\) is an event on \(\Omega\) such that \(\mu(F)=1\) for every such ergodic measure \(\mu\), then \(\mu(F)=1\) for every invariant measure \(\mu\).

One important class of ergodic processes are mixing processes. For a subgroup \(\Gamma\subseteq \Aut(T_d)\) and a function \(f: \Gamma\to \RR\), we write \(\lim_{\gamma\to\infty}f(\gamma)=a\) if for every \(\eeps\gt 0\) and every \(x\in V(T_d)\), there is some \(N\) such that for all \(\gamma\in \Gamma\) satisfying \(d(x,\gamma.x)\ge N\), we have \(|f(\gamma)-a|\lt \eeps\). We say that a \(\Gamma\)-invariant measure \(\mu\) is \textbf{mixing} if for every pair of events \(A,B\), \[
    \lim_{\gamma\to\infty}\mu(A\cap \gamma.B)=\mu(A)\mu(B)
\]
Another example of interest is the following. Let \((\Omega_1,\mathcal{F}_1,\mu_1), (\Omega_2,\mathcal{F}_2,\mu_2)\) be two probability spaces, and \(\Gamma\) a group acting on them. A measurable function \(\phi: \Omega_1\to \Omega_2\) is called \textbf{\(\Gamma\)-equivariant} if it commutes with the group action, i.e. \[
    \phi(\gamma.\omega)=\gamma.\phi(\omega)\qquad \forall \gamma\in \Gamma, \omega\in\Omega_1
\]
If the push-forward measure \(\phi_\ast\mu_1\) coincides with \(\mu_2\), we call it a \textbf{factor} of \(\mu_1\). When \((\Omega_1,\mathcal{F}_1,\mu_1)\) corresponds to an i.i.d. process, we call its factors \textbf{factor of i.i.d.} (FIID) processes. A factor of an ergodic process is itself ergodic; hence FIIDs are ergodic. \\

Another crucial tool for analyzing invariant processes is the celebrated \textbf{Mass-Transport Principle} (MTP):
\begin{theorem} \label{thm:mass-transport-principle}
    With \(\Omega\) as above, let \(F: V\times V\times \Omega\to [0,\infty]\) be a measurable function. We think of \(F(x,y;\omega)\) as an amount of mass sent from \(x\) to \(y\) in the configuration \(\omega\). If \(F\) is \textit{diagonally invariant}, that is, \(F(\phi.x,\phi.y;\phi.\omega)=F(x,y;\omega)\) for all \(\phi\in \Aut(T_d)\), then the expected total mass sent by a vertex equals the expected mass it receives: \[
        \EE\left(\sum_{y\in V}F(o,y;\omega)\right)=\EE\left(\sum_{x\in V}F(x,o;\omega)\right)
    \]
\end{theorem}
 
The MTP was first introduced by Häggström to prove the following remarkable result (see \cite[Theorem 2.4]{benjamini1999group} or\cite{haggstrom1997infinite}):
\begin{theorem} \label{thm:haggstrom-percolation-condition}
    Let \(A\) be an invariant bond percolation on \(T_d\), with marginal probability \(p\) of an edge being open. If \(p\ge \frac{2}{d}\), then the process contains an infinite open component with positive probability.
\end{theorem}
We return to the MTP, in a more general setting, in \Cref{section:random-walk-unimodular-random-networks}.

\subsection{The Howe-Moore Property}

Recall that \(T_d\) is bipartite, let \((P,Q)\) be the partition. Call an automorphism \(\phi\) of \(T_d\) \textbf{parity-preserving} if \(\phi(P)=P\) and \(\phi(Q)=Q\) (equivalently, if \(d_T(\phi(x),x)\) is even for all vertices \(x\)). An automorphism \(\phi\) that's not parity-preserving satisfies \(\phi(P)=Q,\phi(Q)=P\) and is called \textbf{parity-reversing}. Let \(\Aut_+(T_d)\) denote the subgroup of parity-preserving automorphisms of \(T_d\). This is a normal subgroup of \(\Aut(T_d)\), of index \(2\). Note that while \(\Aut(T_d)\)-invariance implies \(\Aut_+(T_d)\)-invariance, \(\Aut_+(T_d)\)-ergodicity is a stronger property than \(\Aut(T_d)\)-ergodicity. However, because the subgroup has index 2, the two properties are intimately linked:

\begin{proposition} \label{prop:equal-mixture}
    Let \(\Omega\) be some space on which \(\Aut(T_d)\) acts, and \(\mu\) an \(\Aut(T_d)\)-ergodic measure.
    \begin{enumerate}
        \item \(\mu\) is an equal mixture of two \(\Aut_+(T_d)\)-ergodic measures \(\mu_1,\mu_2\). Furthermore, given any parity-reversing automorphism \(\tau\), we may take \(\mu_1,\mu_2\) to satisfy \(\mu_2=\tau_\ast\mu_1\).
        \item When \(\Omega=\RR^{E(T_d)}\), \(\mu_1,\mu_2\) have the same marginal distribution as \(\mu\).
    \end{enumerate}
\end{proposition}
\begin{proof}
    If \(\mu\) is already \(\Aut_+(T_d)\)-ergodic, the proposition holds trivially with \(\mu_i=\mu\).
    Otherwise, there exists some \(\Aut_+(T_d)\)-invariant event \(E_1\) such that \mbox{\(0\lt \mu(E_1)\lt 1\)}. Let \(\tau\) be some parity-reversing automorphism. Because \(\Aut_+(T_d)\) is normal, for every \(\phi\in\Aut_+(T_d)\), \[
        \phi(\tau(E_1))=\tau\big((\tau^{-1}\phi\tau)(E_1)\big)=\tau(E_1)
    \]
    hence \(E_2:=\tau(E_1)\) is also \(\Aut_+(T_d)\)-invariant. Since \(\tau^2\in \Aut_+(T_d)\), the events \(E_1\cup E_2,E_1\cap E_2\) are invariant with respect to the full automorphism group \(\Aut(T_d)\). But \[
        \mu(E_1\cup E_2)\ge \mu(E_1)\gt 0,\quad \mu(E_1\cap E_2)\le \mu(E_1)\lt 1
    \]
    and since \(\mu\) is ergodic, the former must have measure \(1\) and the latter measure \(0\). Since \(\mu\) is also \(\tau\)-invariant, we must have \(\mu(E_1)=\mu(E_2)=\frac{1}{2}\). Taking \[
        \mu_i(F):=\mu(F\mid E_i)=2\mu(F\cap E_i)
    \]
    we get \(\mu=\frac{1}{2}\mu_1+\frac{1}{2}\mu_2\). Furthermore \(\mu_2=\tau_\ast\mu_1\), since for any event \(F\), \[
        \tau_\ast\mu_1(F)=2\mu(\tau^{-1}(F)\cap E_1)=2\mu(\tau^{-1}(F\cap E_2))=\mu_2(F)
    \]
    To see that \(\mu_i\) are \(\Aut_+(T_d)\)-ergodic, fix some \(\Aut_+(T_d)\)-invariant event \(F\). By a similar argument to the above, the events \[
        F'=(F\cap E_1)\cup \tau(F\cap E_1),\quad F''=(F\cap E_1)\cap \tau(F\cap E_1)
    \]
    are \(\Aut(T_d)\)-invariant. \(F''\subseteq E_1\cap E_2\) must have measure \(0\), so \[
        \mu(F')=\mu(F\cap E_1)+\mu(\tau(F\cap E_1))=2\mu(F\cap E_1)=\mu_1(F)
    \]
    but \(\mu\) is \(\Aut(T_d)\)-ergodic, hence \(\mu_1(F)=\mu(F')\in \set{0,1}\). The argument for \(\mu_2\) is analogous.

    Now suppose that \(\Omega=\RR^{E(T_d)}\). Since \(\mu_2=\tau_\ast\mu_1\), the distribution of some edge \(e\) under \(\mu_2\) is equal to that of \(\tau(e)\) under \(\mu_1\). But \(\Aut_+(T_d)\) acts transitively on the edges of the tree (since every edge has one endpoint in \(P\) and one in \(Q\)), hence the distribution of \(\tau(e)\) under \(\mu_1\) is equal to that of \(e\) under \(\mu_1\), i.e. the marginal distributions of \(\mu_1,\mu_2\) coincide. Since \(\mu\) is a mixture of the two, its marginal distribution must also coincide with theirs.
\end{proof}

Our interest in \(\Aut_+(T_d)\) stems from a result by Lubotzky and Mozes (\cite{Lubotzky1992}), which says that it satisfies the ``Howe-Moore property'': for every unitary representation of \(\Aut_+(T_d)\), that has no non-zero invariant vectors, its matrix coefficients vanish at infinity (see \cite{ciobotaru2015} for a modern reference). In particular, considering the natural representation of the group on \(L^2(\mu)\) yields the following.

\begin{theorem}[{\cite[Corollary 1]{Lubotzky1992}}] \label{thm:howe-moore-property}
    Let \(\mu\) be a \(\Aut_+(T_d)\)-invariant measure on \(\RR^{E(T_d)}\). If \(\mu\) is \(\Aut_+(T_d)\)-ergodic, then \(\mu\) is mixing.
\end{theorem}

\Cref{thm:howe-moore-property} allows us to justify a claim we made in the introduction: let \(A\) be an \(\Aut(T_d)\)-ergodic process on the edges or vertices of \(T_d\), with marginal distribution \(X\). Fix a ray \(\gamma\), and let \(\phi\in\Aut(T_d)\) act as the shift along \(\gamma\). We claimed that the partial averages of \(A\) along any fixed ray \(\gamma\) converge a.s. to \(\EE(X)\); this follows from the following corollary.

\begin{corollary} \label{corollary:ray-ergodicity}
    Let \(\phi\in \Aut(T_d)\) act as the shift on some ray \(\gamma\). Then any \(\Aut(T_d)\)-ergodic measure \(\mu\) is \(\phi\)-ergodic.
\end{corollary}
\begin{proof}
    Write \(\mu=\frac{1}{2}\mu_1+\frac{1}{2}\mu_2\) with \(\mu_2=\phi_\ast\mu_1\), as in \Cref{prop:equal-mixture} (note that \(\phi\) is indeed parity-reversing). By \Cref{thm:howe-moore-property}, \(\mu_i\) are mixing.
    Let \(E\) be a \(\phi\)-invariant event. Note that \(\phi^2\) is parity-preserving, and that \(d(\phi^{2n}(o),o)=2n\). Therefore, \[
        \mu_i(E)=\mu_i(E\cap \phi^{2n}(E))\xrightarrow[n\to\infty]{}\mu_i(E)^2
    \]
    hence \(\mu_i(E)\in \set{0,1}\). Furthermore, \(\mu_2(E)=\mu_1(\phi^{-1}(E))=\mu_1(E)\), hence \[
        \mu(E)=\frac{1}{2}(\mu_1(E)+\mu_2(E))=\mu_1(E)\in\set{0,1}\qedhere
    \] 
\end{proof}

\subsection{Asymptotics and speeds}
Given some finite subgraph \(K\subseteq T\), we denote by \(A_K\) the sum over all the edges or vertices of \(K\) -- explicitly, \[
    A_K=\sum_{e\in E(K)}A(e)\quad \text{or} \quad A_K=\sum_{v\in V(K)} A(v)
\]
depending on whether \(A\) is an edge or vertex process. In particular we write \(A_{x,y}:=A_{[x,y]}\).
For an edge process \(A\), define the \textbf{lower} and \textbf{upper (asymptotic) means} of a walk \(\mathbf{w}\) as \[
    \MM_\ast(\mathbf{w};A)=\liminf_{n\to\infty}\frac{1}{n}\sum_{i=0}^{n-1} A_{w_{i},w_{i+1}},\quad 
    \MM^\ast(\mathbf{w};A)=\limsup_{n\to\infty}\frac{1}{n}\sum_{i=0}^{n-1} A_{w_{i},w_{i+1}}
\]
When the limit exists, we write \(\MM(\mathbf{w};A)\) for its value. When the process \(A\) is clear from the context, we omit it and write \(\MM_\ast(\mathbf{w}),\MM^\ast(\mathbf{w}),\MM(\mathbf{w})\), respectively. Note that the averages are taken over the entire paths \([w_i, w_{i+1}]\), rather than just the vertices of the walk \((w_i)_{i}\). Similar definitions may be made for vertex processes, but then we have to account for the fact that the end point of each path \([w_{i}, w_{i+1}]\) is counted twice. \\
For a ray \(\gamma\), the definitions above simplify to \[
\MM_\ast(\gamma)=\liminf_{n\to\infty}{\textstyle\frac{1}{n}} A_{\gamma_0,\gamma_n},\quad \MM^\ast(\gamma)=\limsup_{n\to\infty}{\textstyle\frac{1}{n}} A_{\gamma_0,\gamma_n}
\]
In fact, these quantities depend only on the end \(\zeta\) of \(\gamma\): for simplicity, let us assume \(A\) is an edge process. If \(\gamma'\) is equivalent to \(\gamma\), then for all \(n\) \[
    \frac{1}{n}|A_{\gamma_0,\gamma_n}-A_{\gamma'_0,\gamma'_n}|\le \frac{1}{n}|A|_{\gamma\triangle \gamma'}\xrightarrow[n\to\infty]{} 0
\]
where \(|A|(e):=|A(e)|\). Hence the asymptotic means of \(\gamma,\gamma'\) coincide, and the asymptotic means \(\MM_\ast(\zeta),\MM^\ast(\zeta)\) are well-defined. \\

Let us now return to the speeds from \Cref{def:speeds}. In light of the above, \(\sburst\) and \(\ssust\) can be written as \begin{align*}
    \sburst(A)&=\sup_{\zeta\in \partial T}\MM^\ast(\zeta;A) \\ 
    \ssust(A)&=\sup_{\zeta\in \partial T}\MM_\ast(\zeta;A) 
\end{align*}
Presented in this form, it is clearly independent of the choice of root \(o\). The same is true for \(\scloud\):

\begin{claim}
    Let \(o'\in V(T)\). Write \[
    M_n:=\max_{d_T(v,o)=n}A_{o,v},\quad M'_n:=\max_{d_T(v,o')=n}A_{o',v}
\]
Then \(\limsup_n\frac{1}{n}M_n=\limsup_n\frac{1}{n}M'_n\).
\end{claim}
\begin{proof}
Fix a vertex \(v\), and let \(u:=v\wedge_o o'\). Then \(u\) is on both paths \([o,v],[o',v]\), hence \[
    A_{o,v}=A_{o,u}+A_{u,v},\quad A_{o',v}=A_{o',u}+A_{u,v}
\]
Taking the difference, and using the fact that \(u\) is also on the path \([o,o']\), yields \[
    |A_{o,v}-A_{o',v}|=|A_{o,u}-A_{o',u}|\le |A|_{o,o'}
\]
Now for every \(n\), let \(v_n\in S_n(o)\) be a vertex such that \(A_{o,v_n}=M_n\). Write \(r_n=d_T(o',v_n)\). Then \[
\frac{1}{n}M_n\le \frac{1}{n}A_{o',v_n}+\frac{|A|_{o,o'}}{n}\le \frac{r_n}{n}\frac{1}{r_n}M_{r_n}'+o(1)
\]
Now \(\frac{r_n}{n}\to 1\) because \(|r_n-n|=|d_T(o',v_n)-d_T(o,v_n)|\le d_T(o,o')\). Taking \(\limsup\) on both sides of the above, \[
    \limsup_{n\to\infty}\frac{1}{n}M_n
    \le \limsup_{n\to\infty}{\frac{1}{r_n}M'_{r_n}}
\]
By symmetry, equality must hold.
\end{proof}

On \(T_d\), the three speeds are clearly invariant under automorphisms that preserve the root \(o\). Together with the above, we conclude that they are invariant under all automorphisms. In particular for an ergodic process \(A\), they are a.s. constant. \\

Finally, we remark on the measurability of the speeds. Typically, the \(\sigma\)-algebra considered for our probability spaces is the one generated by cylinders; equivalently, it is the Borel \(\sigma\)-algebra induced by the product topology. It is clear that \(\scloud\) is measurable with respect to it. \(\ssust\) is too, however the argument is a bit more subtle: for every \(\alpha\in \RR\), define \[
    V_\alpha:=\set{v\in V(T)\mid A_{o,v}\gt d_T(0,v)\alpha}
\]
and note that \(\ssust\gt \alpha\) if and only if, for some \(\eeps\gt 0\), the graph induced by \(V_{\alpha+\eeps}\) contains an infinite component (which, in particular, contains an infinite ray). Now the events "a vertex \(u\) is in an infinite component of \(V_\alpha\)" are measurable, the set of all vertices is countable and it suffices to consider rational \(\eeps\gt 0\); therefore we have that the sets \(\set{\ssust\gt \alpha}\) are measurable. \\
However, \(\sburst\) is \textit{not} measurable with respect to the Borel \(\sigma\)-algebra. In fact, the sets \(\set{\omega\mid \sburst(\omega)\gt \alpha}\) are \textbf{complete analytic} sets, which are in particular not Borel measurable. That being said, these sets are \textbf{universally measurable}, which means they are measurable with respect to the completion of every Borel probability measure. Thus we will assume that our measurable spaces are complete.  

\section{Equality of \(\scloud\) and \(\sburst\)} \label{section:scloud-equals-sburst}

Given an infinite subtree \(T'\subseteq T\), any process \(A\) on \(T\) restricts to a process on
\(T'\), which we denote \(A|_{T'}\). This restriction can't increase the speed:

\begin{observation}
    Let \(s\) be one of \(\scloud,\sburst,\ssust\). Then \[
    s(A|_{T'})\le s(A)
    \]
\end{observation}

\begin{proof}
    Since the speeds are independent of the choice of root, we may assume that \(T',T\) share the same root \(o\). Then \(\partial T'\subseteq \partial T\) and \(S_r^{(T')}(o)\subseteq S_r^{(T)}(o)\), and the claim follows by definition of the speeds.
\end{proof}

In particular, consider the case where \(T=T_d\) is the \(d\)-regular tree.
Define a \textit{\(d\)-ary tree} as a rooted tree in which each vertex has exactly \(d\) children. \((d-1)\)-ary trees arise naturally as subtrees of \(T_d\) corresponding to directed edges: Let \(\vec{e}=(x;y)\) be a directed edge in \(T_d\). We denote by \(T(\vec{e})=T(x;y)\) the component of \(y\) in \(T\setminus\set{e}\), rooted at \(y\). Then \(T(x;y)\) is a \((d-1)\)-ary tree, which we think of as being directed "away" from \(x\).

\begin{claim}
    For each of the speeds, \[
    s(A)=\max_{x\sim o}s(A|_{T(o;x)})
    \]
\end{claim}
\begin{proof}
    For \(\ssust\), note that \(\partial T=\bigsqcup_{x\sim o}\partial T(o;x)\),
    hence \[
        \ssust(A)=\sup_{\zeta\in \partial T}\MM_{\ast}(\zeta)=\sup_{x\sim o}\sup_{\zeta\in \partial T(o;x)}\MM_\ast(\zeta)=\max_{x\sim o}\ssust(A|_{T(o;x)})
    \]
    Of course, the same argument applies for \(\sburst\). As for \(\scloud\), by the previous claim \(\scloud(A|_{T(o;x)})\le \scloud(A)\) for each \(x\sim o\), so their maximum is also bounded by \(\scloud(A)\). For the other direction, by the definition of \(\scloud\) there a.s. exists a subsequence \((n_k)_k\) such that \[
    \frac{1}{n_k}\max_{|v|=n_k}A_{o,v}\xrightarrow[k\to\infty]{} \scloud(A)
    \]
    as well as a sequence \(v_k\in S_{n_k}(o)\) of vertices that attain the maximal value at each level. By the pigeonhole principle, at least one of the subtrees \(T(o;x)\) contains an infinite number of these vertices, and for that same subtree, \[
    \scloud(A|_{T(o;x)})=\lim_{k\to\infty}\frac{1}{n_k}A_{0,v_k}=\scloud(A)
    \]
\end{proof}

\clearpage

\begin{claim}
    Let \(\vec{e}\) be a directed edge in \(T_d\). Then for each of the speeds, \(s\left(A|_{T(\vec{e})}\right)=s(A)\) a.s.
\end{claim}
\begin{proof}
    In light of the previous claim, it suffices to show that for every \(x_1,x_2\sim o\), \[
        s(A|_{T(o;x_1)})\overset{\as}{=} s(A|_{T(o;x_2)})\textbf{}
    \]
    Since \(d\ge 3\), there's another neighbor \(y\sim o\), different from \(x_1,x_2\). Then \[
        T(o;x_i)\subseteq T(y;o)
    \]
    and hence \(s(A|_{T(y;o)})\overset{\as}{\ge } s(A|_{T(o;x_i)}) \).
    On the other hand, \(s(A|_{T(y;o)})\overset{d}{=} s(A|_{T(o;x_i)})\), and together we conclude that \[
        s(A|_{T(o;x_1)})=s(A|_{T(y;o)})=s(A|_{T(o;x_2)})\qquad \as
    \]
\end{proof}

\begin{corollary} \label{corr:sburst-equals-scloud}
    \(\sburst=\scloud\) a.s.
\end{corollary}

\begin{proof}
    Using ergodic decomposition, it suffices to prove the corollary for an ergodic process \(A\). In this case, \(\scloud(A)=s\) is almost surely constant; by the previous claim, \(\scloud(A|_{T(\vec{e})})\overset{\as}{=}s\) for every directed edge \(\vec{e}\). Since the number of directed edges in \(T_d\) is countable, we also have the equality a.s. for all directed edges at once. \\
    We construct a random ray \(\gamma\) as follows. Start with \(\vec{e}_0=(x_0;y_0)\), where \(y_0=o\) and \(x_0\sim o\) is arbitrary. Given \(\vec{e}_i=(x_i;y_i)\), by the above we have \(\scloud\left(A|_{T\left(\vec{e}_i\right)}\right)=s\); by the definition of the cloud speed, for some distance large enough we have a vertex \(y_{i+1}\in T(\vec{e}_i)\) such that \[
    \displaystyle A_{y_i,y_{i+1}}\ge (s-\textstyle\frac{1}{i})d(y_i,y_{i+1})
    \]
    Define \(\vec{e}_{i+1}:=(x_{i+1};y_{i+1})\), where \(x_{i+1}\) is the predecessor of \(y_{i+1}\) on the path \([y_i, y_{i+1}]\). \\
    By construction, the paths \([y_i, y_{i+1}]\) are edge-disjoint, hence their concatenation \(\gamma\) is an infinite ray. We claim that \(\MM^\ast(\gamma)=s\). Since \((y_i)_i\) is a subsequence of \(\gamma=(\gamma_n)_n\), we have \[
        \MM^\ast(\gamma)\ge \limsup_{i\to\infty}\frac{1}{|y_i|}A_{o,y_i}
    \]
    Fix \(j\in\NN\), then for every \(i\ge j\), by construction \[
        A_{y_{i}, y_{i+1}}\ge (s-\tfrac{1}{i})d(y_i,y_{i+1})\ge (s-\tfrac{1}{j})d(y_i,y_{i+1})
    \] 
    Therefore \[
        \frac{1}{|y_i|}A_{o, y_i}=\frac{1}{|y_i|}A_{o, y_j}+\frac{1}{|y_i|}A_{y_j,y_i}\ge o(1)+\frac{d(y_j,y_i)}{|y_i|}(s-\textstyle\frac{1}{j})\longrightarrow s-\frac{1}{j}
     \]
     Hence \(\MM^\ast(\gamma)\ge s-\frac{1}{j}\); this is true for all \(j\), so \(\MM^\ast(\gamma)\ge s\).
\end{proof}

\section{A Process with \(\scloud\gt \ssust\)} \label{section:worm-process}
We now construct an invariant percolation process for which \(\scloud\gt \ssust\). We define it as an edge process on \(T_3\) for simplicity, for the same process and proof work, \textit{mutatis mutandis}, for a vertex process or on \(T_d\) for \(d\gt 3\).

\subsection{The model}

\begin{definition}[Worms and Territories]
    A \textbf{worm of size \(n\)} is simply a path \(P\) of length \(n\). A \textbf{(worm) territory of size \(n\)} is an \(n\)-buffer \(B_n(P)\) around a worm \(P\) of size \(n\).
\end{definition}

We imagine the buffer \(B_n(P)\), i.e. the territory, as the worm's require "living space". In our setting worms are aggressive toward one another, and if the territories of two worms overlap, they will compete over it. If worms of different size compete, the larger one wins and the smaller one perishes; if both worms are of the same size, they both perish. This motivates the following definition:

Let \(p_n:=c2^{-2n}n^{-10}\), for some small positive constant \(c\). For every vertex \(x\) and every \(n\), toss a coin with probability \(p_n\) of success; if it succeeds, pick a vertex \(y\in S_n(x)\) uniformly, and "spawn" a worm of size \(n\) on the path \([x,y]\) (with the corresponding territory around it). We then denote \(y\) by \(Y_n(x)\) and the territory by \(R_n(x)\) (if a worm is not spawned, let \(R_n(x),Y_n(x):=\emptyset\)). Observe that \(R_n=(R_n(x))_{x\in V}\) is an i.i.d. process, indexed by \(V\) and with values in finite subgraphs of \(T\). We also define the composite process \(R(x):=(R_n(x))_{n=1}^\infty\).

We say a spawned worm \textit{survives} if its territory \(R_n(x)\) does not intersect any other territory \(R_m(x')\) with \(m\ge n\). Define our percolation, \(A\), to be the union of the worms that survive. Note that \(A\) is a factor of the i.i.d. process \(R\), and in particular is automorphism invariant and ergodic.

We also define a "capped" variant of the model, with bounded maximal worm size: denote by \(A^{[m]}\) the model obtained from \((R_n)_{n=1}^m\) by erasing the territories \(R_k(x)\) that intersect another territory of size \(\ell\ge k\), and taking the union of the surviving ones. Note that \(A^{[m]}\) is an \(8\)-dependent process. \\

In the remainder of this subsection, we provide some general lemmas that will be useful for us later on. Throughout the section, we use \(C\) to denote a constant, not necessarily the same at each occurrence, that does not depend on \(n\). \\
First, let \(\Vol_n\) denote the volume of a territory of size \(n\), i.e. \(\Vol_n=|B_n(P)|\) for a path \(P\) of length \(n\). 

\begin{lemma}
    \(\Vol_n=(n+3)2^n-2\le 2n2^n\)
\end{lemma}

\begin{proof}
Write \(P=(x_0,...,x_n)\), and let \(P_k=(x_0,...,x_k)\) be the initial segment of length \(k\) of \(P\). Then \(B_n(P_k)\subseteq B_n(P_{k+1})\) and \[
    B(P_{k+1})\setminus B_n(P_{k})=\set{y\in S_n(x_{k+1})\mid x_k\le x_{k+1}\le y}
\]
In particular, \(|B_n(P_{k+1})|=|B_n(P_{k})|+2^n\); together with the base case \[
|B_n(P_0)|=|B_n(x_0)|=3\!\cdot\!2^n-2
\]
we get \(|B_n(P_k)|=(3+k)2^n-2\). Taking \(k=n\) yields the claim.
\end{proof}

Since smaller worms survive only if they don't intersect larger worms, the following bound will be central to us:

\begin{lemma} \label{lemma:probability-of-being-in-spawned-worm}
    Let \(I_m(x)\) be the event that a vertex \(x\) is in the territory of some spawned worm of size at least \(m\). Then \[
        \PP(I_m(x))\le \sum_{n=m}p_n\Vol_n\le C2^{-m}m^{-9}
    \]
\end{lemma}

\begin{proof}
    We use the mass-transport principle. Define a mass transport \(F\) as follows: every vertex \(x\), for every \(n\ge m\), sends a mass of \(1\) to every \(y\in R_n(x)\). In other words, \[
        F(x,y;\omega)=\sum_{n=m}^\infty \one_{y\in R_n(x)}
    \]

    Then the expected total mass that a vertex \(x\) sends is 
    \begin{align*}
        \EE\left(\sum_{y\in V}F(x,y;\omega)\right)&=\sum_{k=m}^\infty p_k\Vol_k
        \le \sum_{k=m}^\infty c2^{-2k}k^{-10}\cdot 2k2^k \\
        &=C\sum_{k=m}^\infty2^{-k}k^{-9}
        \le C2^{-m}m^{-9}\sum_{j=0}^\infty 2^{-j}=C2^{-m}m^{-9}
    \end{align*}

    By the mass transport principle, this is also the expected mass that a vertex receives. Finally, note that \(I_m(x)=\set{\sum_{z\in V}F(z,x;\omega)\ge 1}\); using Markov's inequality we get \[
        \PP(I_m(x))\le \EE\left(\sum_{z\in V}F(z,x;\omega)\right)\le C2^{-m}m^{-9} \qedhere
    \]
\end{proof}
\vspace{1em}

Next, let \(q_n\) be the probability that some vertex \(x\) spawns a worm of size \(n\) and it survives. It's possible to calculate \(q_n\) exactly, but the following bound is simpler and more elegant:

\begin{proposition} \label{prop:worm-process-probability-of-survival}
    \(\left(1-\frac{C}{n^8}\right)p_n\le q_n\le p_n\)
\end{proposition}

\begin{proof}
    The upper bound is immediate from the definitions. For the lower bound, fix \(x\) and \(n\). Then \begin{equation} \label{eqn:q_n}
       q_n=p_n\cdot\PP(R_n(x)\ \text{survives}\mid R_n(x)\neq \emptyset)  
    \end{equation}

    Note that in terms of the mass transport \(F\) defined in the proof of \Cref{lemma:probability-of-being-in-spawned-worm}, a spawned worm perishes if and only if one of the vertices in its territory \(R_n(x)\neq \emptyset\) receives mass aside from \(R_n(x)\), that is iff
    \[
        \exists {y\in R_n(x)}: \sum_{z\in V}F(z,y;\omega)-\one_{y\in R_n(x)}\gt 0
    \] 
    To estimate this event, define a new process \(R'\) as such: consider the original process \(R\) and an independent copy \(R''\). Define \(R'_{m}(y)=R_m(y)\) for every pair \((y,m)\) \textit{except} for \((x,n)\), for which we define \(R'_n(x)=R''_n(x)\); in other words we start from \(R\), and "forget" the value of \(R_n(x)\), instead replacing it with an independent copy of the same distribution. \\
    Let also \(F'(x,y;\omega)=\sum_{k\ge n}\one_{y\in R'_k(x)}\) be the mass transport corresponding to \(R'\). Then given \(R_n(x)\neq \emptyset\), for every \(y\in R_n(x)\) we have \[
        \sum_{z\in V}F(z,y;\omega)-\one_{y\in R_n(x)}\le \sum_{z\in V}F'(z,y;\omega)
    \]
    because the difference is \(0\) whenever \(y\in R_n'(x)\), and \(1\) otherwise. Therefore, 
    \begin{align*}
             \PP(R_n(x)&\ \text{perishes} \mid R_n(x)\neq \emptyset)
             \le \PP\left(\exists y\in R_n(x):\sum_{z\in V}F'(z,y;\omega)\gt 0\right)\\
             & =\PP\left(\sum_{y\in R_n(x)}\sum_{z\in V}F'(z,y;\omega)\ge 1\right)
             \overset{(1)}{\le} \EE\left(\sum_{y\in R_n(x)}\sum_{z\in V}F'(z,y;\omega)\right) \\
             &\overset{(2)}{=}\Vol_n\cdot\PP(I_n(x))\le \frac{C}{n^8}
    \end{align*}
    where (1) again follows from Markov's inequality and (2) from the fact that \(R\overset{d}{=}R'\) and from \Cref{lemma:probability-of-being-in-spawned-worm}. Putting this in \ref{eqn:q_n} yields \[
        q_n=\big(1-\PP\left(R_n(x)\ \text{perishes} \mid R_n(x)\neq \emptyset\right)\big)p_n\ge \left(1-\frac{C}{n^8}\right)p_n \qedhere
    \] 
\end{proof}

\begin{proposition} \label{prop:p[m]-p-distance-bound}
    let \(p^{[m]}, p\) be the probability that an edge is open in \(A^{[m]}, A\) respectively. Then \[
    \left|p-p^{[m]}\right|\le \frac{C}{m^8}
    \]   
\end{proposition}
\begin{proof}
    Let \(e\) be an arbitrary edge. Then \[
        \left|p-p^{[m]}\right|=\left|\PP(e\in A)-\PP\left(e\in A^{[m]}\right)\right|\le \PP\left(e\in A\triangle A^{[m]}\right)
    \]
    Where \(A\triangle A^{[m]}\) denotes the (pointwise) symmetric difference of the two percolation processes. Now, the event \(\set{e\in A\triangle A^{[m]}}\) is the disjoint union of \(\set{e\in A\setminus A^{[m]}}\) and \(\set{e\in A^{[m]}\setminus A}\). \(e\in A\setminus A^{[m]}\) implies, in particular, that the endpoints \(x,x'\) of \(e\) are contained in some spawned territory of size greater than \(m\); therefore, by \Cref{lemma:probability-of-being-in-spawned-worm}, \[
        \PP(e\in A\setminus A^{[m]})\le \PP(I_m(x))\le C2^{-m}m^{-9}
    \]
    For the second event, given \(e\in A^{[m]}\), let \(\mathcal{R}=R_k(z)\) be the surviving territory that contains \(e\) in \(A^{[m]}\). Then \(e\notin A\) if and only if \(\mathcal{R}\) intersects a spawned territory of size greater than \(m\); since spawned territories of different sizes are independent, and \(\mathcal{R}\) is of size \(k\le m\), we must have \[
        \PP\left(e\notin A\mid e\in A^{[m]}\right)\le \PP\left(\bigcup_{y\in \mathcal{R}}I_{m+1}(y)\right)\le \Vol_n\PP(I_{m+1}(y))\le \frac{C}{m^8}
    \] 
    Together, \[
        \left|p-p^{[m]}\right|\le \PP\left(A\setminus A^{[m]}\right)+\PP\left(e\in A^{[m]}\setminus A\right)\le \frac{C}{m^8} \qedhere
    \]
\end{proof}

\subsection{An upper bound for \(\ssust\)}
We now show that \(\ssust(A)\le \frac{1}{3}\). The idea is that the density of an infinite path comes from the time it spends inside worms; however surviving worms are disjoint, and each contributes a density of at most \(\frac{1}{3}\) to the lower density.

\begin{definition}
    Let \(\mathcal{R}=R_n(x)\) be a territory of size \(n\) around a worm \(W\). A \textbf{crossing} of \(\mathcal{R}\) is a finite path \(\gamma=(\gamma_0,...,\gamma_\ell)\subseteq \mathcal{R}\) such that \(d_T(\gamma_0,W)=d_T(\gamma_\ell,W)=n\), that is \(\gamma\) starts and ends at the "boundary" of the territory.
\end{definition}

For a finite path \(\gamma\) in \(T_3\), define its \textit{density} as the sum along \(\gamma\) divided by its length, i.e. as \(\MM(\gamma):=\frac{1}{|\gamma|}A_\gamma\).

\begin{claim} \label{claim:crossing-density-bound}
    Let \(\gamma=(\gamma_0,...,\gamma_\ell)\) be a crossing of a territory \(\mathcal{R}=R_n(x)\) around a worm \(W\). The density of \(\gamma\) is at most \(\frac{1}{3}\).
\end{claim}

\begin{proof}
    Write \(\gamma_W=\gamma\cap W\) for the segment of the path \(\gamma\) inside the worm \(W\). Then \(A_{\gamma}=A_{\gamma_W}=|\gamma_W|\). If \(\gamma_W=\emptyset\) then \(\MM(\gamma)=0\). Otherwise, \(\gamma_W\) is a path, as an intersection of two paths. Let \(\gamma',\gamma''\) be the segments of \(\gamma\) before reaching \(W\) and after leaving it, respectively. By the definition of a crossing \(d_T(\gamma_0,W)=d_T(\gamma_\ell,W)=n\), hence \(|\gamma'|=|\gamma''|=n\). On the other hand \(|\gamma_W|\le |W|=n\), and together, \[
        \MM(\gamma)=\frac{1}{|\gamma|}A_\gamma=\frac{|\gamma_W|}{|\gamma|}\le \frac{1}{3}
    \]
\end{proof}

\begin{proposition} \label{prop:worm-process-ssust-upper-bound}
    \(\ssust(A)\le \frac{1}{3}\)
\end{proposition}
\begin{proof}
    Let \(\gamma\) be a ray in \(T_3\). Let \(\mathcal{R}_j:=R_{n_j}(x_j)\) be the surviving territories that \(\gamma\) crosses through. If the collection \((\mathcal{R}_j)_j\) is finite, then \((A_{\gamma_0,\gamma_n})_n\) is bounded and hence \(\MM(\gamma)=0\). \\
    Otherwise, let \(a_j,b_j\) be the times where \(\gamma\) enters and leaves each \(\mathcal{R}_j\), respectively. Then \((\gamma_{a_j},...,\gamma_{b_j})\) is a crossing of \(\mathcal{R}_j\). Territories of surviving worms don't overlap, so the intervals \([a_j,b_j]\subseteq \NN\) are pairwise disjoint. Now consider the leaving times \((b_j)_{j=1}^\infty\); since this is a subsequence of \(\NN\), we have \[
        \MM_\ast(\gamma)\le \liminf_{j\to\infty}\frac{1}{b_j}A_{\gamma_0,\gamma_{b_j}}
    \] 
    But the density of \(\gamma\) comes entirely from the time it spends inside territories; together with the previous claim, we obtain \[
        \frac{1}{b_j}A_{\gamma_0,\gamma_{b_j}}=\frac{1}{b_j}\sum_{i=1}^j A_{\gamma_{a_i},\gamma_{b_i}}\le \frac{1}{b_j}\sum_{i=1}^j \frac{1}{3}(b_i-a_i)\le \frac{1}{3}
    \]
\end{proof}

\subsection{A lower bound for \(\scloud\)}

As we've seen in \Cref{claim:crossing-density-bound}, crossing a territory contributes a density of at most \(\frac{1}{3}\). However if we allow ``stopping'' in the middle of the territory, the density can get up to \(\frac{1}{2}\). \\
Our strategy is to look for worms of size \(n\), at distance \(\approx 2n\) from the root \(o\), that face ``away'' from \(o\). Given such a worm, consider the path that starts at \(o\), travels up to the worm, then traverses it. This path (of length \(\approx 3n\)) has \(n\) ones coming from traversing the worm itself, and another \(n\) zeroes from traveling up to the worm inside its territory. However the initial segment of length \(\approx n\), before the path reaches the territory, also contributes some positive, if small, density. Together we get a density a little bigger than \(\frac{1}{3}\).

To describe this formally, we need some notation: \begin{itemize}
    \item Let \(p\) be the marginal probability of an edge being open in \(A\).
    \item Let \(r_n=2n+2\sqrt{n}\) (this is the distance from \(o\) where we'll be looking for worms).
    \item Let \(r_n':=r_n+n=3n+2\sqrt{n}\) (this is the length of the dense paths we'll be looking for, including the worms).
    \item For a vertex \(x\), and \(k\le |x|\), let \(x|_k\) denote the unique vertex on the path \([o,x]\) such that \(d_T(o,x|_k)=k\). In other words, \([o,x|_k]\) is the initial segment of \([o,x]\) of length \(k\).
\end{itemize}

\begin{proposition}\label{prop:worm-process-scloud-lower-bound}
    Let \(\eeps=\frac{p}{60}\). Then a.s. for all but finitely many \(n\), \[
        \frac{1}{r_n'}\max_{|x|=r_n'}A_{o,x}\ge \frac{1}{3}+\eeps
    \]
\end{proposition}

We prove \Cref{prop:worm-process-scloud-lower-bound} via a Borel-Cantelli argument.

\begin{definition}
    Let \(x\in S_{r_n}(o)\). Call \(x\) an \textbf{(\(n\)-)candidate} if \(x\) spawns a worm of size \(n\) that faces away from the root \(o\), i.e. \(R_n(x)\neq \emptyset\) and \(o\le x\le Y_n(x)\). Let \(V_n\subseteq S_{r_n}(o)\) denote the (random) set of candidates, and \(\mu_n:=\EE(|V_n|)\). \\
    
    Also, let \(m_n=\sqrt{n}\). Define the following events: \begin{itemize}
        \item \(D_n\) is the event that \(\frac{1}{r_n'}\max_{|x|=r_n'}A_{o,x}\ge \frac{1}{3}+\eeps\).
        \item \(E_n^1(x)\) is the event that \(x\) is a candidate and the corresponding territory \(R_n(x)\) survives.
        \item \(E_n^2(x)\) is the event that \(x\) is a candidate, and the density of the capped process \(A^{[m_n]}\), along the initial segment \([o,x|_n]\), is at least \(5\eeps\).
        \item \(E_n^3(x)\) is the event that \(x\) is a candidate, and the density of the (pointwise) difference \(A^{[m_n]}\setminus A\), along the initial segment \([o,x|_n]\), is at most \(\eeps\).
    \end{itemize}
    Also, let \[
        E_n(x)=E_n^1(x)\cap E_n^2(x)\cap E_n^3(x), \qquad E_n=\bigcup_{x\in S_{r_n}(o)}E_n(x)
    \]
\end{definition}

\begin{claim} \label{claim:worm-process-E_n-implies-D_n}
    For \(n\) large enough, \(E_n\) implies \(D_n\).
\end{claim}

\begin{proof}
    Suppose \(x\in S_{r_n}(o)\) satisfies \(E_n^1(x)\cap E_n^2(x)\cap E_n^3(x)\), and let \(y:= Y_n(x)\). Then \([o,y]\) is a path of length \(r_n+n=r_n'\). Since \(E_n^1(x)\) occurs, the worm \([x,y]\) survives and contributes a density of \(1\). In addition, because \(E_n^2(x)\) and \(E_n^3(x)\) occur, the initial segment \([o,x|_n]\) contributes a density of \[
        \frac{1}{n} A_{o,x|_n}\ge \frac{1}{n}\big(A^{[m_n]}\big)_{o,x|_n}-\frac{1}{n}\big(A^{[m_n]}\setminus A\big)_{o,x|_n}\ge 4\eeps
    \]
    Therefore, \[
        \frac{1}{r_n'}A_{o,y}\ge \frac{1}{r_n'}(n\cdot 4\eeps +n\cdot 1)=\frac{1+4\eeps}{3+o(1)}\xrightarrow[n\to\infty]{} \frac{1+4\eeps}{3}\gt \frac{1}{3}+\eeps
    \]
    Hence for large enough \(n\), \(D_n\) is implied.
\end{proof}

\begin{claim} \label{claim:worm-process-E_n-bound}
    Let \(x\in S_{r_n}(o)\) be an arbitrary vertex. Then \begin{equation} \label{eqn:worm-process-E_n-bound}
        \PP(E_n^c)\le \PP\left(|V_n| \lt \frac{\mu_n}{2}\right)+2\sum_{i=1}^3 \PP\left(\left(E_n^i(x)\right)^c\mid x\in V_n\right)        
    \end{equation}
\end{claim}

\begin{proof}
     Let \[
        U_n^i:=\set{x\in V_n\mid E_n^i(x)\text{ does not occur}}
     \]
     and \(U_n:=U_n^1\cup U_n^2\cup U_n^3\). Then \(E_n^c=\set{U_n=V_n}\), hence \[
        \PP(E_n^c)=\PP(|U_n|=|V_n|)\le \PP\left(|V_n|\lt \frac{\mu_n}{2}\right)+\PP\left(|U_n|\ge \frac{\mu_n}{2}\right)
     \]
     For the latter, by Markov's inequality, \[
        \PP\left(|U_n|\ge \frac{\mu_n}{2}\right)\le 2\frac{\EE(|U_n|)}{\mu_n}\le \frac{2}{\mu_n}\sum_{i=1}^3\EE(|U_n^i|)
     \]
     Now \(\mu_n=\EE(|V_n|)=|S_{r_n}(o)|\PP(x\in V_n)\), where \(x\in S_{r_n}(o)\) is arbitrary. Similarly, \begin{align*}
        \EE(|U_n^i|)&=|S_{r_n}(o)|\PP(x\in U_n^i) \\
        &=|S_{r_n}(o)|\PP(x\in V_n)\cdot \PP(x\in U_n^i\mid x\in V_n) \\
        &=\mu_n\PP(x\in U_n^i\mid x\in V_n)=\mu_n\PP\left(\left(E_n^i(x)\right)^c\;\middle |\;  x\in V_n\right)
     \end{align*}
     Therefore \begin{align*}
         \PP(E_n^c)&\le \PP\left(|V_n|\lt \frac{\mu_n}{2}\right)+\frac{2}{\mu_n}\sum_{i=1}^3\EE(|U_n^i|) \\
         &\le \PP\left(|V_n|\lt \frac{\mu_n}{2}\right)+2\sum_{i=1}^3\PP\left(\left(E_n^i(x)\right)^c\;\middle |\;  x\in V_n\right)\qedhere
     \end{align*}
\end{proof}

We proceed by bounding the probabilities on the right side of \eqref{eqn:worm-process-E_n-bound}.

\begin{claim} \label{claim:worm-process-V_n-bound} 
\(\PP\left(|V_n|\lt \frac{\mu_n}{2}\right)\le C2^{-\sqrt{n}}\).
\end{claim}

\begin{proof}
    Observe that \(|V_n|\sim \operatorname{Bin}(|S_{r_n}(o)|, \frac{2}{3}p_n)\). As a result, its expected value is \[
        \mu_n=|S_{r_n}(o)|\cdot\frac{2}{3}p_n=C2^{2\sqrt{n}}n^{-10}
    \] 
    and \(\Var(|V_n|)\le \mu_n\). By Chebychev's inequality, \[
        \PP\left(|V_n|\lt \frac{\mu_n}{2}\right)\le 4\frac{\Var(|V_n|)}{\mu_n^2}\le \frac{4}{\mu_n}=\frac{Cn^{10}}{2^{2\sqrt{n}}}\le \frac{C}{2^{\sqrt n}} \qedhere
    \]
\end{proof}
\vspace{1em}

\begin{claim} \label{claim:worm-process-E_n^1-bound}
    \(\PP\left(\left(E_n^1(x)\right)^c\;\middle |\;  x\in V_n\right)\le \frac{C}{n^8}\).
\end{claim}
\begin{proof}
    Note that given \(x\in V_n\), \(\left(E_n^1(x)\right)^c\) occurs if and only if the spawned territory \(R_n(x)\) fails to survive; applying \Cref{prop:worm-process-probability-of-survival} yields the given bound.
\end{proof}

\begin{claim} \label{claim:worm-process-E_n^2-bound}
    For \(n\) large enough, \[
        \PP\left(\left(E_n^2(x)\right)^c\;\middle |\; x\in V_n\right)
        =\PP\left(\big(A^{[m_n]}\big)_{o,x|_n}\lt 5\eeps n\right)\le e^{-C\sqrt{n}}
    \]
\end{claim}

\begin{proof}
    For the first equality, observe that \(\left(E_n^2(x)\right)^c\) is the intersection of the two events \(\set{x\in V_n}\) and \(\left\{\big(A^{[m_n]}\big)_{o,x|_n}\lt 5\eeps n\right\}\); these are independent, since the former depends only on \(R_n\), and the latter only on \((R_j)_{j=1}^{m_n}\).

    Let \(p^{[m]}, p\) be as in \Cref{prop:p[m]-p-distance-bound}, then for \(m\) large enough, \(p^{[m]}\ge \frac{p}{2}\). Recall that \(A^{[m]}\) is an \(8m\)-dependent process and \(\eeps=\frac{p}{60}\).

    Let \(\ell_n=\floor{10m_n}=\floor{10\sqrt{n}}\), and write \[
        n=b_n\ell_n+\rho, \quad 0\le \rho \lt \ell_n
    \]
    Divide the initial segment \([o,x|_n]\) into \(b_n=\floor{\frac{n}{\ell_n}}\) blocks of length \(\ell_n\), plus a remainder of length \(\rho\) which we ignore. Because \(A^{[m_n]}\) is \(8m_n\)-dependent and \(\ell_n\gt 8m_n\), any two blocks which are separated by at least one other block are independent; thus we restrict our attention to even-indexed blocks.
    Let \(k=\floor{\frac{b_n}{2}}\) be the number of even blocks. For \(1\le j\le k\), let \(B_j\) be the number of open edges in the \(2j\)-th block, and \(Y_j=\frac{1}{\ell_n}B_j\) be its density. Then \((Y_j)_j\) are independent and identically distributed, with \(Y_j\in [0,1]\) and \(\EE(Y_j)=p^{[m_n]}\ge \frac{p}{2}\).

    Applying Hoeffding's inequality to \((Y_j)_j\), we obtain \begin{equation} 
        \begin{split}
            \PP\left(\frac{1}{k}\sum_{j=1}^k Y_j\le \frac{p}{4} \right)&\le \PP\left(\frac{1}{k}\sum_{j=1}^k Y_j-\EE(Y_j)\le -\frac{p}{4}\right) \\
            &\le e^{-2k\left(\tfrac{p}{4}\right)^2}\le e^{-C\sqrt{n}}
        \end{split}
        \label{eqn:worm-process-E_n^2-bound-1}
    \end{equation}
    where the last inequality follows from the fact that \(k=\Theta(\sqrt{n})\).

    Now, given that the even blocks do achieve a density of at least \(\frac{p}{4}\), the density of the initial segment \([o,x|_n]\) is bounded from below by \[
        \frac{1}{n}\big(A^{[m_n]}\big)_{o,x|_n}\ge \frac{1}{n}\sum_{j=1}^k B_j=\frac{\ell_n}{n}\sum_{j=1}^k Y_j\ge \frac{\ell_nk}{n}\frac{p}{4}
    \]
    We can bound \(\ell_n k\) from below: \[
        \ell_nk=\ell_n\floor*{\frac{b_n}{2}}\ge \ell_n\left(\frac{b_n}{2}-1\right)\ge \frac{n}{2}-\left(\frac{\rho}{2}+\ell_n\right)=\frac{n}{2}-o(n)
    \]
    Thus for large enough \(n\), \begin{equation} \label{eqn:worm-process-E_n^2-bound-2}
        \frac{1}{n}\big(A^{[m_n]}\big)_{o,x|_n}\ge \frac{\left(\tfrac{n}{3}\right)}{n}\frac{p}{4}\ge \frac{p}{12}=5\eeps
    \end{equation}

    Combining \eqref{eqn:worm-process-E_n^2-bound-1} and \eqref{eqn:worm-process-E_n^2-bound-2}, we conclude \[
        \PP\left(\big(A^{[m_n]}\big)_{o,x|_n}\lt 5\eeps n\right)\le \PP\left(\frac{1}{k}\sum_{j=1}^k Y_j\le \frac{p}{4} \right)\le e^{-C\sqrt{n}} \qedhere
    \]  
\end{proof}
\vspace{1em}

\begin{claim} \label{claim:worm-process-E_n^3-bound}
    \(\PP\left(\left(E_n^3(x)\right)^c\;\middle |\;  x\in V_n\right)\le \frac{C}{n^4}\).
\end{claim}
\begin{proof}
    Fix \(x\in S_{r_n}(o)\), and write \([o,x|_n]=(e_1,e_2,...,e_n)\). By Markov's inequality, we have \begin{equation}
        \label{eqn:worm-process-E_n^3-bound-1}
        \begin{aligned}
            \PP\left(\left(E_n^3(x)\right)^c\;\middle |\; x\in V_n\right)
            &=\PP\left(\left(A^{[m_n]}\setminus A\right)_{o,x|_n}\ge \eeps n\;\middle |\; x\in V_n\right) \\
            &\le \frac{\EE\left(\left(A^{[m_n]}\setminus A\right)_{o,x|_n}\;\middle |\; x\in V_n\right)}{\eeps n} \\
            &=\frac{1}{\eeps n}\sum_{j=1}^n\PP\left(e_j\in A^{[m_n]}\setminus A\;\middle |\; x\in V_n\right)
        \end{aligned}
    \end{equation}

    Fix some \(1\le j\le n\), and write \(e:=e_j\) for brevity. We show that \begin{equation}
        \label{eqn:worm-process-E_n^3-bound-2}
        \PP(e\in A^{[m_n]}\setminus A\mid x\in V_n)\le \PP(e\in A^{[m_n]}\setminus A)
    \end{equation}
    To see this, write the left side as \begin{align*}
        \PP(e\in A^{[m_n]}\setminus A\mid x\in V_n)&=\PP(e\in A^{[m_n]}, e\notin A\mid x\in V_n) \\
        &=\PP(e\notin A\mid e\in A^{[m_n]}, x\in V_n)\PP(e\in A^{[m_n]}\mid x\in V_n) \\
        &=\PP(e\notin A\mid e\in A^{[m_n]}, x\in V_n)\PP(e\in A^{[m_n]})
    \end{align*}
    where the last equality follows since the events \(\set{e\in A^{[m_n]}}, \set{x\in V_n}\) are independent. \\
    Now given \(e\in A^{[m_n]}\), let \(\mathcal{R}=R_k(z)\) be the territory that \(e\) is part of. Then \(e\notin A\) if and only if some spawned territory of size greater than \(m_n\) intersects \(\mathcal{R}\). Conditioning on \(x\in V_n\) means that the territory \(R_n(x)\) cannot intersect \(\mathcal{R}\) (since \(R_n(x)\) faces away from the root); spawned territories are independent, hence \[
        \PP(e\notin A\mid e\in A^{[m_n]}, x\in V_n)\le \PP(e\notin A\mid e\in A^{[m_n]})
    \]
    Multiplying both sides by \(\PP(e\in A^{[m_n]})\) yields \eqref{eqn:worm-process-E_n^3-bound-2}. \\
    Plugging \eqref{eqn:worm-process-E_n^3-bound-2} into \eqref{eqn:worm-process-E_n^3-bound-1}, and applying \Cref{prop:p[m]-p-distance-bound}, we obtain \begin{align*}
        \PP\left(\left(E_n^3(x)\right)^c\;\middle |\; x\in V_n\right)
        &\le \frac{1}{\eeps}\PP(e\in A^{[m_n]}\setminus A) \\
        &\le\frac{1}{\eeps}\PP(e\in A^{[m_n]}\triangle A)\le \frac{C}{m_n^8}=\frac{C}{n^4} \qedhere
    \end{align*}
\end{proof}
\vspace{1em}

\begin{proof}[Proof of \Cref{prop:worm-process-scloud-lower-bound}]
    We combine the previous results. \Cref{claim:worm-process-E_n-bound} and the following four bounds together imply that \[
        \PP(E_n^c)\le \PP\left(|V_n| \lt \frac{\mu_n}{2}\right)+2\sum_{i=1}^3 \PP\left(\left(E_n^i(x)\right)^c\mid x\in V_n\right) \le \frac{C}{n^4}
    \]
    Now by \Cref{claim:worm-process-E_n-implies-D_n}, we conclude that \[
        \sum_{n=1}^\infty \PP(D_n^c)\le \sum_{n=1}^\infty \PP(E_n^c)\lt \infty
    \]
    By the Borel-Cantelli lemma, a.s. all but finitely many of the events \((D_n)_n\) must occur.
\end{proof}

With the combination of Propositions \ref{prop:worm-process-ssust-upper-bound} and \ref{prop:worm-process-scloud-lower-bound}, we've established that \(\scloud\) and \(\ssust\) do not coincide. However, given the equality between \(\scloud\) and \(\sburst\) that we've established in \Cref{corr:sburst-equals-scloud}, one might ask whether there's some variation on \(\scloud\) that does coincide with \(\ssust\). The immediate candidate that comes to mind is replacing the limsup in \(\scloud\) by a liminf -- in other words, define \[
     \underline{\scloud}(A):=\liminf_{n\to\infty}\frac{1}{n}\max_{|x|=n}A_{o,x}
\]
Does the equality \(\underline{\scloud}(A)=\ssust(A)\) hold? It turns out that our process \(A\) provides a negative counterexample to that as well.

\begin{proposition} \label{prop:worm-process-scloud-variant-bound}
    For the process \(A\) discussed above, and \(\eeps\) as in \Cref{prop:worm-process-scloud-lower-bound}, we have \[
        \underline{\scloud}(A)\ge \frac{1}{3}+\eeps
    \]
    In particular, \(\scloud(A)\neq \ssust(A)\).
\end{proposition}

In fact, \Cref{prop:worm-process-scloud-variant-bound} is a simple consequence of \Cref{prop:worm-process-scloud-lower-bound} and the following lemma (which we state without proof):

\begin{lemma} \label{lemma:subsequence-bdd-diff}
    Let \((a_n)_n\) be a sequence with bounded differences, i.e. \(|a_{n+1}-a_n|\le C\) for some \(C\) independent of \(n\). 
    Let \((a_{n_k})_k\) be a subsequence such that \((n_k)_k\) has bounded differences. Then \[
        \liminf_{n\to\infty}\frac{a_n}{n}=\liminf_{n\to\infty}\frac{a_{n_k}}{n_k}
    \] 
\end{lemma}

\begin{proof}[Proof of \Cref{prop:worm-process-scloud-variant-bound}]
    Let \(a_n:=\max_{|x|=n} A_{o,x}\). Observe that \(a_{n+1}\le a_n+1\): indeed, let \(z\in S_{n+1}(o)\) be a vertex for which \(A_{o, z}=a_{n+1}\), and let \(y\in S_n(o)\) be the preceding vertex on the path \([o, z]\). Then \[
    a_{n+1}=A_{o,z}=A_{o,y}+A_{y,z}\le \max_{x\in S_n(o)}A_{o,x}+1=a_n+1
    \]
    Together with the non-negativity of \(A\), this implies \(|a_{n+1}-a_n|\le 1\). 
    Additionally, let \(n_k=r_k'\), then \[
        r_{k+1}'-r_k=3+(2\sqrt{k+1}-2\sqrt{k})\le 5
    \]
    \Cref{lemma:subsequence-bdd-diff} now implies \[
        \underline{\scloud}=\liminf_{n\to\infty}\frac{a_n}{n}=\liminf_{n\to\infty}\frac{a_{n_k}}{n_k}\ge \frac{1}{3}+\eeps
    \]
    where the inequality follows from \Cref{prop:worm-process-scloud-lower-bound}.
\end{proof}

\section{Unimodular Networks and Random Walks on Trees} \label{section:random-walk-unimodular-random-networks}

In this section we develop a few tools that will be needed in the next section, but are also interesting in their own right. The main results in this section are \Cref{corollary:rwrn-speed} and \Cref{thm:mean-along-limiting-end}, which identify the speed of random walks on \textit{unimodular random networks} (URNs) supported on trees, and the asymptotic means along their limiting ends. The framework presented here is based primarily on the work of Aldous and Lyons (\cite{aldous2018processes}). \\

\subsection{Preliminaries}

We start with some definitions. A \textbf{network} is a locally finite, possibly directed graph \(G=(V,E)\) with a map \(c: E\to [-\infty,\infty]\). We call the values of this map \textbf{marks} or \textbf{weights}. A \textbf{rooted network} \((G,o)\) is a network \(G\) with a distinguished vertex \(o\). An \textbf{isomorphism} of rooted networks \((G_1,o_1),(G_2,o_2)\) is an isomorphism of the underlying graphs \(\varphi: G_1\to G_2\) such that \(\varphi(o_1)=o_2\). We generally do not distinguish between a rooted network and its isomorphism class. \\
Let \(\gstar\) denote the space of all (isomorphism classes of) connected rooted networks. This space separable and complete when endowed with the following metric. We say \((G_1,o_1),(G_2,o_2)\in \gstar\) are \textit{\(r\)-close}, for some \(r\gt 0\), if there exists an isomorphism \(\phi\) of rooted graphs between the balls of radius \(\floor{r}\) around the roots \(o_i\) in \(G_i\), such that for every pair or corresponding edges or vertices, the distance between their marks is less than \(1/r\). Define the distance of \((G_1,o_1),(G_2,o_2)\) to be \(1/(1+\alpha)\), where \[
    \alpha=\sup\set{r\ge 0: (G_1,o_1),(G_2,o_2)\text{ are \(r\)-close}}
\]
Similarly, we define the space \(\gstarr\) of all (equivalence classes of) connected networks with two distinguished vertices, with a topology similar to \(\gstar\).

\begin{definition}
    Let \(\mu\) be a measure on \(\gstar\). We call \(\mu\) \textbf{unimodular} if it satisfies the \textbf{Mass-Transport Principle}: for any Borel function \(f: \gstarr\to [0,\infty]\), \[
        \int_{\gstar} \sum_{y\in V(G)}f(G,o,y)\ d\mu([G,o])=\int_{\gstar} \sum_{x\in V(G)}f(G,x,o)\ d\mu([G,o])
    \]
    A \textbf{unimodular random network} is a random element of \(\gstar\) whose law is a unimodular measure. 
\end{definition}

To see the connection between unimodular measures and automorphism-invariant processes, let \(G\) be a fixed Cayley graph of a finitely-generated group. Then any \(\Aut(G)\)-invariant measure \(\nu\) corresponds naturally to a unimodular measure \(\mu\), by choosing a fixed root \(o\). In particular, the law of an automorphism-invariant process on \(G\) can be represented as a unimodular measure. A more general correspondence, beyond the transitive case, can be found in \cite[Theorem 3.2]{aldous2018processes}. \\

Next, we discuss random walks. Let \(G\) be some fixed network. Broadly speaking, a \textbf{random walk} is a random element of \(V(G)^\NN\). However, we will only be interested in random walks that arise as Markov chains, with \(V(G)\) as their state space. Recall that a \textbf{Markov chain} with (countable) state space \(S\) is a random \(S\)-valued process \((X_n)_{n=0}^\infty\) such that the distribution of \(X_{n+1}\), given the past, depends only on the value of \(X_n\), or more explicitly \[
    (X_{n+1}\mid X_1=x_1,...,X_n=x_n)\overset{d}{=}(X_{n+1}\mid X_n=x_n)
\]
for all \(n\). We call \(p(x,y):=\PP(X_{n+1}=y\mid X_n=x)\) the \textbf{transition probabilities} of the chain (we consider only homogeneous chain, in which the transition probabilities do no depend on \(n\)). As an example, consider \textbf{simple random walk} on \(G\), whose transition probabilities are given by \[
    p_G(x,y)=\frac{1}{\deg_G(x)}\one_{y\sim x}
\]
A \textbf{(\(p\)-)stationary measure} for a Markov chain is a function \(\pi: S\to [0,\infty]\) such that for all \(x\in S\), \begin{equation} \label{eqn:stationarity-condition}
    \pi(y)=\sum_{x\in S}\pi(x)p(x,y)   
\end{equation}
\(\pi\) is called \textbf{(\(p\)-)reversible} with respect to the chain if for all \(x,y\in S\), \begin{equation} \label{eqn:reversibility-condition}
    \pi(x)p(x,y)=\pi(y)p(y,x)  
\end{equation}
Note that \(p\)-reversibility implies \(p\)-stationarity, and that both are preserved when \(\pi\) is scaled by a positive constant. Reversible Markov chains may be identified with a particular type of random walks: Let \(G\) be a network with positive marks. Define \[
\pi(x)=\sum_{y\in V(G)}c(x,y)
\]   
and consider the Markov chain whose transition probabilities are proportional to the marks, that is \[
    p(x,y)=\frac{c(x,y)}{\pi(x)}
\]
We call this chain the \textbf{network random walk} corresponding to \(G\). Note that it is in particular reversible. Conversely, given a Markov chain with state space \(S\) and reversible distribution \(\pi\), define a network whose vertices are \(S\), and where two vertices \(x,y\) are connected iff \(p(x,y)\gt 0\), with mark \[
c(x,y)=\pi(x)p(x,y)=\pi(y)p(y,x)
\]
Then \(p(x,y)=c(x,y)/\pi(x)\), and the chain is given the structure of a network random walk. We call the underlying graph of this network the \textbf{supporting graph} of the chain. \\

Given a chain with a stationary measure \(\pi\), we define the \textbf{reversed Markov chain} by the transition probabilities \[
    \hat{p}(x,y):=\frac{\pi(y)}{\pi(x)}p(y,x)
\]
As its name suggests, the reversed Markov chain is obtained from the original chain by going backward in time instead of forward; more precisely, the distribution of \(X_i\) given the ``future'' \((X_{i+1},X_{i+2},...)\), is given by the reversed transition probabilities \(\hat{p}(x,y)\). In particular, since the reversed process is also a Markov chain, it follows that \[
    (X_i\mid X_{i+1}=x_{i+1},X_{i+2}=x_{i+2},...)\overset{d}{=}(X_{i}\mid X_{i+1}=x_{i+1})
\]
The reversed chain allows us to extend the original chain backwards: given \(X_0\), define \(X_{-1}\) as a step of the reversed chain starting at \(X_0\). Then \((X_{-1},X_0,X_1,...)\) is also a Markov chain given by the transition probabilities \(p\). \\
Finally, note that \(\hat{p}=p\) if and only if the Markov chain is \(\pi\)-reversible. This implies, in particular, that for reversible chains, \[
    (X_{-1}\mid X_0=x_0)\overset{d}{=}(X_1\mid X_0=x_0)
\]
We emphasize that reversibility is not required for the reversed chain to be well-defined; only a stationary measure is required. \\

We also consider random walks on random networks. Here a network \((G,o)\) is first drawn from some unimodular measure \(\mu\) on \(\gstar\), then a random walk is run on \(G\), starting from \(o\). Of course, in this case the random walk may interact with the underlying environment in different ways. In a \textbf{random walk in random environment} (RWRE), the network's (random) marks determine the walk's transition probabilities. More precisely, a \textbf{random environment} is a Borel function \(p: \gstarr\to [0,1]\), written as \(p: (G,x,y)\mapsto p_G(x,y)\), such that \(\sum_{y\in V}p_G(x,y)=1\) for all \(x\in G\). Alternatively, the walk may observe the marks around its position as it traverses the space, in which case we call the marks a \textbf{random scenery}. The two may also be combined into a random walk in a random environment with random scenery (RWRERS).

Since our model now has two stochastic ``layers'' -- the network and the walk on it -- we distinguish between two types of measures: for each network \((G,o)\) in the support of \(\mu\), let \(\theta_{(G,o)}\) denote the distribution of the random walk on \((G,o)\); in other words, these distributions treat the environment as deterministic, and only the walk is random. We call these the \textbf{quenched} distributions. The \textbf{annealed} distribution, on the other hand, is obtained by integrating the quenched distributions over all possible networks.
More explicitly, let \(\pstar\) be the space of (equivalence classes of) pairs \(\big((G,x_0),(x_n)_{n=0}^\infty\big)\) with \((G,x_0)\in \gstar\) and \(x_n\in V(G)\), with a topology similar to \(\gstar\) and \(\gstarr\). Then the annealed measure \(\Theta\) is given by \begin{equation} \label{eqn:annealed-measure}
    \Theta(F)=\int_{\pstar}\ \theta_{(G,o)}(F)\ d\mu(G,o)
\end{equation}
for all events \(F\). To illustrate the difference, suppose we put i.i.d. marks \(c(x,y)\sim U([0,1])\) on the edges of the \(3\)-regular tree \(T_3\). Define a random environment by \[
    p(x,y)=\begin{cases}
        1 & y\sim x,\  c(x,y)=\max_{z\sim x} c(x,z) \\ 
        0 & \text{otherwise}
    \end{cases}
\]
That is, the corresponding walk always picks the edge, adjacent to its current location, with the maximal mark. Note that since the distribution of the marks is continuous, the maximum is a.s. unique for all \(x\), hence \(p\) is well-defined as a random environment. Then the quenched behavior of the walk, i.e. its distribution on a fixed realization of the marks, is completely deterministic; in particular two walks on the same environment will always take their first step over the same edge. However in the annealed case the distribution of the first step is uniform over the neighbors of the starting point \(x_0\), since the distribution of the marks is invariant under automorphisms. \\

To be able to analyze the walk \((x_n)_{n}\) and the underlying network, we would like the \textit{environment viewed from the particle}, i.e. the sequence \(\set{(G,x_n)}_n\), to be stationary under the \textbf{shift} \(\sigma: \pstar\to\pstar\), defined by \[
    \sigma\big((G,x_0),(x_n)_n\big):=\big((G,x_1),(x_{n+1})_{n}\big)
\]
This seldom happens automatically, even in simple cases:
for example, let \(G\) be some finite graph (which we identify with a network with \(c\equiv 1\)). Construct a random rooted network by choosing the root \(o\) uniformly across the vertices of \(G\), then run simple random walk \((x_n)_n\) starting from that root. Then \[
    \frac{1}{n}\sum_{k=0}^{n-1} \deg_G (x_n)\xrightarrow[n\to\infty]{\as} \frac{\EE(\deg_G(o)^2)}{\EE(\deg_G(o))}
\]
(this follows from \Cref{thm:rwrers} below). If the sequence were stationary, we would expect the sequence of averages on the left to converge to \(\EE(\deg_G o)\); but by Jensen's inequality, the term on the right hand side is at least \(\EE(\deg_G o)\), with equality if and only if the \(\deg_G (o)\) is constant, that is iff \(G\) is regular. In other words, over time the walk tends to land on vertices with higher degrees more often than others. 

However, when the quenched walks have a stationary distribution, we can construct a new measure with respect to which the walk is stationary by \textbf{biasing} with respect to it: let \(\mu\) be unimodular, and \(p\) be some random environment defined on the support of \(\mu\). Also let \(\nu: \gstar\to (0,\infty)\) be a Borel map, written as \(\nu: (G,x)\to \nu_G(x)\). We say that \(\nu\) is \textbf{\(p\)-stationary} if for all \(G\), \(\nu_G\) is a stationary distribution for the random walk on \(G\) determined by \(p_G\). \textbf{\(p\)-reversibility} is defined analogously. Define a measure \(\Theta_{\mu,p,\nu}\) on \(\pstar\) by \[
    \Theta_{\mu,p,\nu}(F):=\int_{\pstar}\theta_{(G,o)}(F)\nu_G(o)\ d\mu(G,o)
\]
That is, \(\Theta_{\mu,p,\nu}\) is the annealed distribution of the random walk determined by \(p_G\), starting at \(o\) and biased by \(\nu_G(o)\), where \((G,o)\) is distributed according to \(\mu\). In this context, we call \(\nu\) the \textbf{bias}. Note that \[
    \Theta_{\mu,p,v}(\pstar)=\int_{\gstar} \nu_G(o)\ d\mu(G,o)
\]
hence if \(\int \nu\ d\mu=1\), \(\Theta_{\mu,p,\nu}\) is a probability measure; more generally if \(c:=\int \nu\ d\mu\lt\infty\), then the normalized bias \(\tilde\nu:=\frac{1}{c}\nu\) yields a probability measure. The following theorem establishes stationarity with respect to the biased measure (see \cite[Theorem 4.1]{aldous2018processes}, \cite{lyons1999stationary} for the proof):

\begin{theorem} \label{thm:rwrers}
    Let \(\mu, p\) and  \(\nu\) be as above. If \(\nu\) is \(p\)-stationary, then the biased annealed measure \(\Theta_{\mu,p,\nu}\) is stationary under the shift \(\sigma\). If \(\nu\) is \(p\)-reversible, then \(\Theta_{\mu,p,\nu}\) is reversible, that is for all events \(F_0,F_1\), \begin{equation}
        \Theta_{\mu,p,\nu}\big[(G,x_0)\in F_0,(G,x_1)\in F_1\big]
        =\Theta_{\mu,p,\nu}\big[(G,x_1)\in F_0,(G,x_0)\in F_1\big]   
    \end{equation}
\end{theorem}

Note that \(\Theta_{\mu,p,\nu}\) is absolutely continuous with respect to the annealed measure \(\Theta\) from (\ref{eqn:annealed-measure}), with the bias \(\nu\) as the Radon-Nikodym derivative. Since \(\nu\) is strictly positive, the two measures are in fact equivalent. In particular, \(\Theta(F)=1\) for some event \(F\) if and only if \(\Theta_{\mu,p,\nu}(F)=1\), and thus we may state that an event occurs a.s. without specifying under which of the two. \\

We call this model -- of a RWRERS on a unimodular random network, with a stationary bias \(\nu\) -- a \textbf{random walk on a random network} (RWRN). We emphasize that that this is more general than a \textit{network random walk} on a unimodular random network; the former is in particular reversible, while a RWRN need not be. As an example, let \(G\) be any unimodular random graph (that is, a unimodular random network with all marks equal to \(1\) a.s.). Consider \textit{non-backtracking random walk} on \(G\): this is the Markov chain whose state space is directed edges of \(G\), and with transition probabilities \[
    p_G((x;y),(y;z))=\frac{1}{\deg_G(y)-1}\one_{z\neq x}
\]
This walk is stationary with respect to the bias \(\nu\equiv 1\). However, it is non-reversible. \\

Next, we present a sufficient condition for ergodicity of a RWRN under the shift \(\sigma\). For such a condition to be general, we would expect it to include assumptions on both the environment and the walk itself. First, observe that the set of unimodular measures on \(\gstar\) is a convex set; we call a unimodular measure \(\mu\) \textbf{extremal}, or \textbf{ergodic}, if it cannot be written as a non-trivial convex combination of other unimodular measures. Also, let \(\mathcal{I}\) denote the \(\sigma\)-algebra of events on \(\gstar\) that are invariant under non-rooted isomorphisms. 

\begin{theorem}[{\cite[Theorem 4.7]{aldous2018processes}}]
    Let \(\mu\) be a unimodular measure on \(\gstar\). Then \(\mu\) is extremal if and only if it is \(\mathcal{I}\)-trivial, i.e. \(\mu(I)\in\set{0,1}\) for all \(I\in\mathcal{I}\).
\end{theorem}

\begin{theorem}[{\cite[Theorem 4.6]{aldous2018processes}}] \label{thm:rwrn-ergodicity-criterion}
    Let \(\mu\) be an extremal unimodular measure on \(\gstar\), \(p\) a random environment, and \(\nu\) a \(p\)-stationary bias. Suppose that \(\mu\)-a.s. for every \(G\) and for every pair of neighbors \(x\sim y\) in \(G\), we have \(p_G(x,y)\gt 0\). Then the biased measure \(\Theta_{\mu,p,\nu}\) is ergodic under the shift.
    
\end{theorem}

Finally, for the sake of brevity we introduce ``probabilistic'' variants of the notation defined above for quenched, annealed and biased measures. Given some random variable \(X\), its expected value is taken, by default, with respect to the \textit{annealed} measure \(\Theta\): \[
    \EE(X):=\int X\ d\Theta
\]
The underlying network \(G\) will typically be denoted as \(\omega\); the quenched expected value of \(X\) on \(\omega\) is then denoted \(\EE(X\mid \omega)\). \\
Given a bias \(\nu\), the expected value of \(X\) with respect to the \textit{normalized} biased measure \(\Theta_{\mu,p,\tilde{\nu}}\) is denoted \(\EE_\nu(X)\). Alternatively, \[
    \EE_\nu(X):=\EE(X\tilde\nu(o))=\frac{\EE(X\nu(o))}{\EE(\nu(o))}
\]
Similarly, for an event \(F\), \(\PP(F)\) denotes its probability with respect to the annealed measure, and \[
\PP(F\mid \omega):=\EE(\one_F\mid \omega),\quad\PP_\nu(F):=\EE_\nu(\one_F)
\]

\subsection{Asymptotic means along the limiting end}
Let \(\xx\) be an ergodic, end-convergent RWRN, supported on an ergodic random tree \((T,o)\) with marks \(A\). We assume throughout that the sum of \(A\) along a single step has finite expected value, i.e. \(\EE_\nu(|A_{x_0,x_1}|)\lt \infty\). Define \[
    \Delta_{i,j}:=A_{x_i,x_j},\qquad i\le j
\]
This process is stationary, in the sense that its distribution is invariant under the shift \(\sigma: \Delta_{i,j}\mapsto \Delta_{i+1,j+1}\) (with respect to the biased measure); our assumptions also imply that \(\Delta\) is ergodic and \(\EE_{\nu}(|\Delta_{0,1}|)\lt\infty\).
If \(A\ge 0\), then the process \(\Delta\) is \textit{subadditive}, meaning that \(\Delta_{i,k}\le \Delta_{i,j}+\Delta_{j,k}\) for all \(i\le j\le k\). In this case, by Kingman's subadditive ergodic theorem, the limit \[
    \tau_A:=\lim_{n\to\infty}\frac{\Delta_{0,n}}{n}
\]
exists a.s. and in \(L^1\) (see e.g. \cite[Theorem 14.44]{lyonsperes}). Furthermore, since \(\Delta\) is ergodic, \(\tau_A\) is a.s. constant and \[
    \tau_A=\lim_{n\to\infty}\frac{1}{n}\EE_\nu(\Delta_{0,n})
\]
If \(A\) is not nonnegative, \(\Delta\) might not be subadditive; nevertheless the same results of the subadditive ergodic theorem follow by writing \[
    \Delta=\Delta^+-\Delta^-,\qquad \Delta_{i,j}^\pm:=(A^\pm)_{x_i,x_j}
\]
and noting that \(\Delta^\pm\) are ergodic, subadditive and \(\EE_\nu(\Delta_{0,1}^\pm)\le \EE_\nu(|\Delta_{0,1}|)\lt\infty\). \\

\(\tau_A\) is called the \textbf{time constant} of the process \(\Delta\). An important particular example is \(A\equiv 1\), where \(\Delta_{i,j}=d(x_i,x_j)\) and the time constant is \[
    \tau_1=\lim_{n\to\infty}\frac{d(x_0,x_n)}{n}=\speed(\xx)
\] 

In general, the time constant can be difficult to estimate. One approach, which was in fact introduced by Kingman in his original proof of the subadditive theorem (\cite{kingman1968}), is to find a \textit{subadditive decomposition} for \(X\), that is a pair of stationary processes \(Y,Z\) such that \begin{itemize}
    \item \(X=Y+Z\)
    \item \(Y\) is additive with the same time constant as \(X\). \\
    (Note that since \(Y\) is additive, its time constant is simply \(\EE(Y_{0,1})\))
    \item \(Z\) is nonnegative and subadditive with time constant \(0\).
\end{itemize} 

Kingman proved that such a decomposition always exists (it is generally not unique). As we've mentioned in our setting \(\Delta\) might not be subadditive, but the essence of the subadditive decomposition remains -- in fact, the fact that the walk is supported on trees gives it an explicit, elegant form:

\begin{theorem} \label{thm:busemann-function}
Let \(\xx\to\xi\) be an end-convergent, ergodic RWRN on a unimodular random tree \((T,A,o)\). The limit \begin{equation} 
    \label{eqn:busemann-function-definition}
        B_{i,j}:=\lim_{N\to\infty}(\Delta_{i,N}-\Delta_{j,N}),\qquad i\le j 
    \end{equation}
    a.s. exists, and equals \begin{equation} \label{eqn:busemann-function-form}
        B_{i,j}=A_{x_i,x_j}-2A_{x_j,x_i\wedge_{x_j} \xi}\qquad \as
    \end{equation}
    Furthermore \(B\) is stationary, additive and has the same time constant as \(\Delta\), that is, \begin{equation} \label{eqn:busemann-function-time-constant}
        \EE_\nu(B_{0,1})=\tau_A
    \end{equation}
\end{theorem}

\begin{remark}
   If \(\Delta\) is subadditive, then by the triangle inequality, \[
        B_{i,j}=\lim_{N\to\infty}(\Delta_{i,N}-\Delta_{j,N})\le \Delta_{i,j}
    \]
    Hence, \(\Delta-B\) is nonnegative and subadditive with time constant \(0\), such that \(\Delta=B+(\Delta-B)\) is a proper subadditive decomposition. 
\end{remark}

\begin{proof}[Proof of \Cref{thm:busemann-function}]
    Note that for all \(i,j,N\), \[
        A_{x_i,x_j}+A_{x_j,x_N}=A_{x_i,x_N}+2A_{x_j,x_i\wedge_{x_j} x_N}
    \]
    Rearranging gives \[
        \Delta_{i,N}-\Delta_{j,N}=A_{x_i,x_j}-2A_{x_j,x_i\wedge_{x_j} x_N}
    \]
    We claim, then, that \(A_{x_i,x_j}-2A_{x_j,x_i\wedge_{x_j} x_N}\to A_{x_i,x_j}-2A_{x_j,x_i\wedge_{x_j} \xi}\). We show that more generally, for vertices \(x,y\) and a walk \(z_n\to \zeta\), \[
        A_{y,x\wedge_y z_n}\xrightarrow[n\to\infty]{} A_{y,x\wedge_y \zeta}
    \]
    Indeed, let \(r:=d(x,y)\) and \(\gamma=[y,\zeta)\). By \ref{eqn:end-convergence-geometric-condition}, for \(n\) large enough we have \[
        y\le \gamma_r\le z_n
    \]
    This implies that for large enough \(n\), \(x\wedge_y z_n=x\wedge _y\gamma_r\); this is because the path \([x,y]\) is of length \(r\), so \[
        [y,x]\cap [y,z_n]=[y,x]\cap [y,\gamma_r]
    \]
    and in particular the vertex farthest from \(y\) is the same in both intersections. By the same argument, \(x\wedge_y \zeta=x\wedge_y \gamma_r\), so for all \(n\) large enough, \[
        x\wedge_y z_n=x\wedge_y \zeta
    \]
    and in particular \(A_{y,x\wedge_y z_n}\to A_{y,x\wedge_y \zeta}\).
    
    Thus, the limit \eqref{eqn:busemann-function-definition} a.s. exists, and has the form \eqref{eqn:busemann-function-form}. Next we prove its properties. For stationarity, we have \[
        B_{i+1,j+1}=\lim_{N\to\infty}(\Delta_{i+1,N}-\Delta_{j+1,N})\overset{d}{=}\lim_{N\to\infty}(\Delta_{i,N-1}-\Delta_{j,N-1})=B_{ij}
    \]
    Where the second equality follows from the fact that \(\Delta\) is stationary under the shift, and that if \(X_N\overset{d}{\to} X, Y_N\overset{d}{\to }Y\) are sequences of random variables such that \(X_N\overset{d}{=}Y_N\) for all \(N\), then \(X\overset{d}{=}Y\). \\

    For additivity, let \(i\le j\le k\), then \begin{align*}
        B_{i,j}+B_{j,k}&=\lim_{N\to\infty}(\Delta_{i,N}-\Delta_{j,N})+\lim_{N\to\infty}(\Delta_{j,N}-\Delta_{k,N}) \\
        &=\lim_{N\to\infty}(\Delta_{i,N}-\Delta_{k,N})=B_{i,k}
    \end{align*}
    Finally, write \(\EE_\nu(\Delta_{0,n})=\sum_{k=0}^{n-1}\EE_\nu(\Delta_{0,k+1}-\Delta_{0,k})\) as a telescopic sum, then \[
        \tau_A=\lim_{n\to\infty}\frac{1}{n}\EE_\nu(\Delta_{0,n})=\lim_{n\to\infty}\frac{1}{n}\sum_{k=1}^{n}\EE_\nu(\Delta_{0,k}-\Delta_{0,k-1})
    \]
    On the other hand, \[
        \EE_\nu(B_{0,1})=\lim_{k\to\infty}\EE_\nu(\Delta_{0,k}-\Delta_{1,k})=\lim_{k\to\infty}\EE_\nu(\Delta_{0,k}-\Delta_{0,k-1})
    \]
    Where the first equality follows from dominated convergence, and the second from stationarity of \(\Delta\). \eqref{eqn:busemann-function-time-constant} now follows from the fact that the Cesàro sums of a convergent sequence converge to the same limit.
\end{proof}

\begin{remark}
    The assumption of the walk being supported on a tree is used to guarantee the existence of the limit \eqref{eqn:busemann-function-definition} (and its explicit form), which is the main obstacle for general subadditive processes (see \cite[Section 2.6]{50yearsfpp}). As evident from the proof, once the limit is known to exist, the properties of \(B\) follow from general considerations. \(B\) is called the \textbf{generalized Busemann function} for the process \(\Delta\).
\end{remark}

An immediate corollary of \Cref{thm:busemann-function} is the following formula for the speed of a RWRN on a unimodular tree.

\begin{corollary} \label{corollary:rwrn-speed}
    Let \(\xx\to \xi\) be an end-convergent, ergodic RWRN on a unimodular random tree \((T,o)\). Then \begin{equation} \label{eqn:rwrn-speed}
        \speed(\xx)=\EE_\nu(d(x_0,x_1))-2\EE_\nu(|x_{-1}\wedge \xi|)\qquad \as
    \end{equation}
\end{corollary}

\begin{proof}
    Applying \Cref{thm:busemann-function} to the process \(A\equiv 1\) and using stationarity, \begin{align*}
        \speed(\xx)&=\EE_\nu(B_{0,1})=\EE_\nu(B_{-1,0}) \\
        &=\EE_\nu(d(x_{-1},x_0))-2\EE_\nu(|x_{-1}\wedge \xi|) \\
        &=\EE_\nu(d(x_0,x_1))-2\EE_\nu(|x_{-1}\wedge \xi|) \qedhere  
    \end{align*}
\end{proof}

\begin{remark}
    \Cref{corollary:rwrn-speed} generalizes Proposition 4.9 in \cite{aldous2018processes}, which says that in the case of simple random walk \(\xx\) on an infinite, ergodic, unimodular random tree \((T,o)\) with \(\EE(\deg_T(o))\lt \infty\), \begin{equation} \label{eqn:aldous2018-speed-formula}
        \speed(\xx)=1-\frac{2}{\EE(\deg_T(o))}\qquad\as    
    \end{equation}

    Indeed, the bias for this walk is given by \(\nu_T(x):=\deg_T(x)\). In addition, \mbox{\(|x_{-1}\wedge \xi|=\one_{\xi_{1}=x_{-1}}\)} and \(d(x_0,x_1)\equiv1\).
    Putting this into \eqref{eqn:rwrn-speed}, we get \[
    \speed(\xx)=1-2\frac{\EE(\one_{\xi_1=x_{-1}}\cdot \deg_To)}{\EE(\deg_To)}\qquad \as
    \]
    But given \(\deg_T(o)=d\), \(x_{-1}\) is independent of \(\xi_1\), hence they are equal with probability \(\frac{1}{d}\). Therefore \[
    \EE(\one_{\xi_1=x_{-1}}\cdot \deg_To)=\EE\left(\frac{1}{\deg_T o}\deg_To\right)=1
    \]
    which yields \eqref{eqn:aldous2018-speed-formula}.
\end{remark}

\Cref{thm:busemann-function} can further be used to derive a formula for the asymptotic mean of \(A\) along the limiting end \(\xi\), i.e. for \(\MM(A;\xi)\) (and in particular, to show that it a.s. exists as a limit):

\begin{theorem} \label{thm:mean-along-limiting-end}
Let \(\xx\to \xi\) be an end-convergent, ergodic RWRN on a unimodular random tree \((T,A,o)\), with positive speed. Then a.s. \begin{equation} \label{eqn:mean-along-limiting-end}
    \MM(\xi;A)=\frac{\EE_\nu(A_{o,x_1})-2\EE_\nu(A_{o,x_{-1}\wedge \xi})}{\speed(\xx)}
\end{equation}
\end{theorem}
\begin{proof}
    In light of \Cref{thm:busemann-function}, the numerator on the right hand side of \eqref{eqn:mean-along-limiting-end} is simply the time constant \(\tau_A\). Let \[
        y_n:=x_n\wedge \xi=o\wedge_{x_n}\xi
    \]
    By \eqref{eqn:busemann-function-form}, a.s. \[
        B_{0,n}=A_{x_0,x_n}-2A_{x_n,y_n}
    \]
    But \begin{align*}
        &\frac{1}{n}A_{x_0,x_n}=\frac{1}{n}\Delta_{0,n}\xrightarrow[n\to\infty]{\as} \tau_A \\
        &\frac{1}{n}B_{0,n}=\frac{1}{n}\sum_{i=0}^{n-1}B_{i,i+1}\xrightarrow[n\to\infty]{\as} \EE_\nu(B_{0,1})=\tau_A
    \end{align*}
    where the first limit is the definition of \(\tau_A\), and the second follows by the ergodic theorem and \Cref{thm:busemann-function}. Therefore \[
        \frac{1}{n}A_{x_n,y_n}\xrightarrow[n\to\infty]{\as} 0
    \]
    On the other hand we have \[
        A_{o,x_n}=A_{o,y_n}+A_{y_n,x_n}
    \]
    Hence \[
        \lim_{n\to\infty}\frac{1}{n}A_{o,y_n}=\lim_{n\to\infty}\frac{1}{n}A_{o,x_n}=\tau_A
    \] \\
    
    Now, note that the supporting tree \((T,\one, o)\), obtained from \((T,A,o)\) by ``forgetting'' the marks \(A\), is also unimodular. In particular, repeating the argument for \(A\equiv 1\), we obtain that \[
        \speed(\xx)=\lim_{n\to\infty}\frac{|y_n|}{n}=\lim_{n\to\infty}\frac{|x_n\wedge \xi|}{n}
    \]
    Let \(\gamma=[o,\xi)\) be the walk's escape ray, then \(y_n=\gamma_{k_n}\) for \(k_n=|y_n|=|x_n\wedge \xi|\). With the above, we get that \[
        \tau_A=\lim_{n\to\infty}\frac{k_n}{n}\frac{1}{k_n}A_{o,\gamma_{k_n}}=\speed(\xx)\lim_{n\to\infty}\frac{1}{k_n}A_{o,\gamma_{k_n}}
    \]
    Hence, we're left with showing that \[
        \MM(\xi;A)=\lim_{k\to\infty}\frac{1}{k}A_{o,\gamma_k}=\lim_{n\to\infty}\frac{1}{k_n}A_{o,\gamma_{k_n}}
    \]
    This would be immediate, expect we don't know apriori that \(\MM(\xi;A)\) exists as a limit. We deduce it from the existence of the limit on the right. First, note that by the ergodic theorem, \[
        \frac{k_{n+1}}{k_n}=\frac{n+1}{n}\frac{k_{n+1}}{n+1}\frac{n}{k_n}\xrightarrow[n\to\infty]{\as}1
    \]
    For every \(k\), choose \(n=n(k)\) be such that \[
        k_n\le k\le k_{n+1}
    \]
    The choice of \(n(k)\) might not be unique, and it might not be strictly increasing. However \(n(k)\to \infty\) as \(k\to\infty\), hence \[
        1\le \frac{k}{k_{n(k)}}\le \frac{k_{n(k)+1}}{k_{n(k)}}\xrightarrow[k\to\infty]{\as} 1
    \]
    Now, write \(S_k:=A_{0,\gamma_k}\), then \begin{align*}
        \left|\frac{1}{k}S_k-\frac{1}{k_n}S_{k_n}\right|&= \left|\frac{1}{k}S_k-\frac{1}{k}S_{k_n}+\frac{1}{k}S_{k_n}-\frac{1}{k_n}S_{k_n}\right| \\
        &\le \frac{1}{k}|S_k-S_{k_n}|+\left|\frac{k_n}{k}-1\right|\cdot\frac{1}{k_n}S_{k_n}
    \end{align*}
    The second term goes to \(0\) since \(\frac{k}{k_n}\to 1\) and \(\frac{1}{k_n}S_{k_n}\to \frac{\tau_A}{\speed(\xx)}=:L\) is bounded. For the first term, write \begin{align*}
        |S_k-S_{k_n}|&=|A_{o,\gamma_k}-A_{o,\gamma_{k_n}}|=|A_{\gamma_{k_n},\gamma_k}| \\
        &\le |A|_{\gamma_{k_n},\gamma_{k_{n+1}}}=|A|_{o,\gamma_{k_{n+1}}}-|A|_{o,\gamma_{k_n}}
    \end{align*}
    Hence \[
        \frac{1}{k}|S_k-S_{k_n}|\le \frac{1}{k_n}\left(|A|_{o,\gamma_{k_{n+1}}}-|A|_{o,\gamma_{k_n}}\right)\xrightarrow[k\to\infty]{} \frac{1}{\speed(\xx)}(\tau_{|A|}-\tau_{|A|})=0
    \]
    Therefore, \[
        \left|\frac{1}{k}S_k-L\right|\le \left|\frac{1}{k}S_k-\frac{1}{k_n}S_{k_n}\right|+\left|\frac{1}{k_n}S_{k_n}-L\right|\longrightarrow 0
    \]
    and thus \(\MM(\xi;A)=\lim\frac{1}{k}S_k=L=\frac{\tau_A}{\speed(\xx)}\), as claimed.
\end{proof}

\begingroup
\renewcommand{\clearpage}{} 
\section{Existence of Exceptional Rays} \label{section:existence-of-exceptional-rays}

In this section we prove Theorems \ref{thm:exceptional-rays-exist} and \ref{thm:vertex-exceptional-rays-exist}, establishing the existence of exceptional rays for both edge and vertex processes. The edge case is obtained by constructing a suitable RWRN, and analyzing it using the tools from \Cref{section:random-walk-unimodular-random-networks}, whereas the vertex case follows as a corollary from the edge case. 

Before addressing \Cref{thm:exceptional-rays-exist} in its full generality, we prove a weakened version via an elegant and relatively simple argument which utilizes Häggström's condition for invariant percolation (\Cref{thm:haggstrom-percolation-condition}) and the central limit theorem. This approach fails in the general case, but it contains a few important elements that will make another appearance there.

\begin{proposition} \label{prop:exceptional-rays-haggstrom-argument}
    Let \(A\) be an ergodic edge process on \(T_d\), \(d\ge 5\). Suppose its marginal distribution \(X\) is non-constant and has finite variance. Then a.s. the burst speed satisfies \[
        \sburst(A)\ge \EE(X)+\delta_X
    \]
    where \(\delta_X\gt 0\) depends only on the law of \(X\).
\end{proposition}

\begin{proof}
    Without loss of generality, we assume that \(\EE(X)=0\). Let \((X_i)_{i}\) be independent random variables distributed as \(X\), and \(Y_k=\frac{1}{k}\sum_{i=1}^k X_i\). By the central limit theorem, for all \(\theta\in \RR\), \[
        \PP\left(Y_k\ge \frac{\theta\sigma_X}{\sqrt{k}}\right)\xrightarrow[k\to\infty]{} \PP(Z\ge \theta)
    \] 
    where \(Z\sim N(0,1)\) and \(\sigma_X\) is the standard deviation of \(X\). Choose \(\theta\gt 0\) such that \(\PP(Z\ge \theta)\gt\frac{2}{5}\), and \(k\) large enough such that \(\PP\left(Y_k\ge \frac{\theta\sigma_X}{\sqrt{k}}\right)\ge \frac{2}{5}\). \\
    Now let \(A_1,...,A_k\) be independent copies of \(A\), and \[
        B(e):=\frac{1}{k}(A_1(e)+...+A_k(e))
    \]
    be their (pointwise) average. Let \(a:=\frac{\theta\sigma_X}{\sqrt{k}}\gt 0\) serve as a ``threshold'', and define an edge percolation \(C\) by \[
        C(e):=\one_{B(e)\ge a}
    \]
    Note that \(B\) and \(C\) are invariant processes, and the marginal distribution of \(B\) is \(Y_k\). Hence the marginal probability of an edge being open in \(C\) is \[
        \PP(B(e)\ge a)=\PP(Y_k\ge a)\ge \frac{2}{5}\ge \frac{2}{d}
    \]
    By Häggström's theorem (\Cref{thm:haggstrom-percolation-condition}), \(C\) contains an infinite component with positive probability; such a component contains an infinite ray \(\gamma\), and by definition of \(C\) every edge \(e\in \gamma\) satisfies \(B(e)\ge a\). In particular \(\sburst(B)\ge a\) with positive probability (note that \(B\) is not necessarily ergodic, despite \(A\) being ergodic).
    Finally, observe that \(\sburst(A)\ge \sburst(B)\) a.s., since for every ray \(\gamma\), \begin{align*}
        \MM^\ast(\gamma;B)&=\limsup_{n\to\infty}\left(\frac{1}{n}B_{\gamma_0,\gamma_n}\right)
        =
        \frac{1}{k}\left(\limsup_{n\to\infty}\frac{1}{n}\sum_{i=1}^k(A_i)_{\gamma_0,\gamma_n}\right) \\
        &\le \frac{1}{k}\sum_{i=1}^k \limsup_{n\to\infty}\left(\frac{1}{n}(A_i)_{\gamma_0,\gamma_n}\right)=\frac{1}{k}\sum_{i=1}^k \MM^\ast(\gamma;A_i) \\
        &\le \frac{1}{k}\sum_{i=1}^k \sburst(A_i)=\sburst(A)
    \end{align*}
    where the last equality follows since \(A\) is ergodic, hence \(\sburst(A), \sburst(A_i)\) are all a.s. constant. Thus with positive probability, \[
        \sburst(A)\ge \sburst(B)\ge a
    \]
    But \(A\) is ergodic, hence \(\sburst(A)\ge a\) a.s. 
\end{proof}

The two differences between \Cref{prop:exceptional-rays-haggstrom-argument} and the full \Cref{thm:exceptional-rays-exist} -- the assumption that \(d\ge 5\) and the result involving \(\sburst\) rather than \(\ssust\) -- are harder to overcome using Häggström's theorem. Regarding the former, note that for \(d=3,4\) the threshold in Häggström's theorem is \(\frac{2}{3}\) or \(\frac{1}{2}\), respectively, and in particular there's no positive \(\theta\) for which \(\PP(Z\ge \theta)\ge \frac{2}{d}\). As for the latter, note that our argument gives, in fact, that \(\ssust(B)\ge a\) (for \(B\) and \(a\) as in the proof); however, attempting to ``translate'' this back to \(A\) fails, since the analogous calculation yields \(\ssust(A)\le \ssust(B)\).

We solve both of these issues by replacing Häggström's theorem with a suitable RWRN, which we investigate in the next subsection. However, our general strategy will still involve deriving another invariant process \(B\) as the average of i.i.d. copies of \(A\), analyzing it, then translating the results back to \(A\). \\

\subsection{The eager random walk}
Let \(n\in \NN\). Define a graph \((G_n)_n\) as follows. \begin{itemize}
    \item For odd \(n\), define \(G_n\) as the graph whose vertices are \(V(G_n)=V(T_d)\), and where \(x\sim y\) if and only if \(d_{T}(x,y)=n\). 
    \item For even \(n\), the same construction gives a graph with two isomorphic connected components; define \(G_n\) as the component that contains the root \(o\).
\end{itemize}
Alternatively, recall that \(T_d\) is the Cayley graph of the group \(\langle a_1,...,a_d\mid a_i^2=e\rangle\) with respect to the generators \(a_1,...,a_d\). Then \(G_n\) is the Cayley graph generated by the set of all reduced words of length \(n\). In particular \(G_n\) is a connected, transitive graph of degree \(d_n:=|S_n(o)|=d\cdot(d-1)^{n-1}\). \\

These graphs appear naturally as the supporting graphs of RWRNs on \(T_d\) whose step size is \(n\). In particular, simple random walk on \(G_n\) is equivalent to simple random walk on \(T_d\) with step size \(n\), i.e. the Markov chain
where \mbox{\((X_{n+1}\mid X_n=x)\)} is distributed uniformly on \(S_n(x)\); the latter is clearly transient, hence \(G_n\) is a transient graph for all \(n\). 

\begin{definition}
    Let \(n\in \NN\), \(\eeps\gt 0\) and \(C\) be an invariant edge percolation on \(G_n\). The \textbf{eager random walk} \(\xx=\xx^{(n,\eeps,C)}\) is defined as follows. Fix a realization \(\omega\) of \(C\). Given that \(x_t=x\), toss a coin with probability \(\frac{\deg_C(x)}{\deg_C(x)+\eeps}\) of success. If the toss succeeds, choose \(x_{t+1}\) uniformly among the neighbors of \(x\) in \(C\);  otherwise, choose \(x_{t+1}\) uniformly among the neighbors of \(x\) in \(G_n\). We call the former an \textbf{eager} step, and the latter a \textbf{mixing} step. For simplicity we take the walk to start at \(x_0=o\).
\end{definition}

The parameter \(n\) is the walk's step size. The parameter \(\eeps\) and the mixing step ensure that the walk is transient and ergodic (see below). The percolation \(C\) encodes a condition on edges of \(G_n\), or equivalently, a condition on paths of length \(n\) in \(T_d\); paths satisfying the condition are heavily favored by the walk. \\

We can model this walk as a RWRN on \(T_d\). Place weights on paths of the tree by \begin{align*}
    c_\omega(x,y)&=\one_{(x,y)\in C}+\frac{\eeps}{d_n}\quad \text{for}\quad d_T(x,y)=n \\
    \nu_\omega(x)&=\sum_{y\in S_n(x)}c_\omega(x,y)=\deg_C(x)+\eeps
\end{align*}

Then the eager random walk \(\xx=\xx^{(n,\eeps,C)}\) is given by the transition probabilities \[
    p_\omega(x,y)=\frac{c_\omega(x,y)}{\nu_\omega(x)}
\]
That is, \(\xx\) has the form of a network random walk on a URN. In particular, it is reversible with respect to \(\nu\). Some additional properties are given in the next observation.

\begin{observation} \label{observation:eager-rw-basic-properties} \leavevmode
    \begin{itemize}
        \item \(p_\omega(x,y)\ge \frac{\eeps}{d_n}\gt 0\) for all \(d_T(x,y)=n\); the graphs \(G_n\) are transient, so it follows from Rayleigh's monotonicity principle that \(\xx\) is transient (see \cite[p. 52]{lyonsperes}).
        \item The walk has bounded step size, so transience is equivalent to end-convergence; thus \(\xi=\xi(\xx)\) is well defined a.s.
        \item If the percolation \(C\) is ergodic, then the walk \(\xx\) is also ergodic, by \Cref{thm:rwrn-ergodicity-criterion}. \\
    \end{itemize}
\end{observation}

Suppose now that the walk \(\xx\) is ergodic. Applying \Cref{thm:mean-along-limiting-end}, we get that for an ergodic process \(B\), a.s. \[
    \MM(\xi;B)=\frac{\EE_\nu(B_{x_0,x_1})-2\EE_\nu(B_{o,x_{-1}\wedge \xi})}{\speed(\xx)}
\]
as long as \(\speed(\xx)\gt 0\). In the remainder of this subsection, we develop an upper bound for the ``penalty term'' \(\EE_\nu(B_{o,x_{-1}\wedge \xi})\) and a lower bound for the speed. Let \(p\) be the marginal probability of \(C\) (i.e. the probability that some fixed edge is open). Note that \(\EE(\nu(o))=d_np+\eeps\).

\begin{lemma} \label{lemma:eager-rw-backtrack-exponential-decay}
For every \(1\le \ell\le n\), we have \[
    \PP_\nu(|x_{-1}\wedge \xi|\ge \ell)\le \frac{1}{p}\frac{1}{(d-1)^\ell}
\]\end{lemma}

\begin{proof}
    Observe that \(|x_{-1}\wedge \xi|\ge \ell\) if and only if \(x_{-1},\xi\) both continue the same vertex \(v\in S_\ell(o)\), that is, iff \(v\) is on both paths \([o,x_{-1}]\) and \([o,\xi)\) (we write these below as \(v\le x_{-1}\), \(v\le \xi\) respectively). \(x_{-1}\) and \(\xi\) are independent given \(\omega\), hence \[
        \PP(|x_{-1}\wedge \xi|\ge \ell\mid \omega)=\sum_{v\in S_\ell(o)}\PP(v\le x_{-1}\mid \omega)\PP(v\le \xi\mid \omega)
    \]
    Therefore \begin{equation}
        \label{eqn:penalty-term-eqn}
        \nu_\omega(o)\PP(|x_{-1}\wedge \xi|\ge \ell\mid\omega)=\sum_{v\in S_\ell(o)}\left(\sum_{v\le x\in S_n(o)}c_\omega(o,x)\right)\PP(v\le \xi\mid \omega)
    \end{equation}
    Now, \[
        \sum_{v\le x\in S_n(o)}c_\omega(o,x)\le \left(1+\frac{\eeps}{d_n}\right)\#\set{x\in S_n(o)\mid v\le x}
    \]
    and routine calculation yields that \[
        \#\set{x\in S_n(o)\mid v\le x}=(d-1)^{n-\ell}=\frac{d_n}{d_\ell}
    \]
    Hence, \[
        \sum_{v\le x\in S_n(o)}c_\omega(o,x)\le\left(1+\frac{\eeps}{d_n}\right)\frac{d_n}{d_\ell}=\frac{d_n+\eeps}{d_\ell}
    \]
    and putting this into \eqref{eqn:penalty-term-eqn}, \[
        \nu_\omega(o)\PP(|x_{-1}\wedge \xi|\ge \ell\mid\omega)\le \frac{d_n+\eeps}{d_\ell}\sum_{v\in S_\ell(o)}\PP(v\le \xi\mid \omega)=\frac{d_n+\eeps}{d_\ell}
    \]
    Finally, taking expectations, \begin{align*}
        \PP_\nu(|x_{-1}\wedge\xi|\ge \ell)&=\frac{\EE(\nu_\omega(o)\PP(|x_{-1}\wedge\xi|\ge \ell\mid \omega))}{\EE(\nu(o))} \\
        &\le \frac{d_n+\eeps}{d_np+\eeps}\frac{1}{d_\ell}\le\frac{1}{p}\frac{1}{d_\ell}\le \frac{1}{p}\frac{1}{(d-1)^\ell}
    \end{align*}
\end{proof}

An immediate application of \Cref{lemma:eager-rw-backtrack-exponential-decay} is the following lower bound on the speed of \(\xx\).

\begin{corollary} \label{corollary:eager-rw-speed-lower-bound}
\[
    \EE_\nu(|x_{-1}\wedge \xi|)\le \frac{1}{p}
\]
As a result, \[
    \speed(\xx)\ge n-\frac{2}{p}
\]
\end{corollary}

\begin{proof}
    Using the tail-sum formula for expectation and \Cref{lemma:eager-rw-backtrack-exponential-decay}, \begin{align*}
        \EE_\nu(|x_{-1}\wedge \xi|)
        &=\sum_{\ell=1}^n\PP_\nu(|x_{-1}\wedge \xi|\ge \ell)
        \le \frac{1}{p}\sum_{\ell=1}^n \frac{1}{(d-1)^{\ell}} \\
        &\le \frac{1}{p}\sum_{\ell=1}^\infty \frac{1}{(d-1)^{\ell}} = \frac{1}{p}\frac{1}{d-2}
        \le \frac{1}{p}     
    \end{align*}
    The second inequality follows from \Cref{corollary:rwrn-speed}.
\end{proof}

\begin{proposition} \label{prop:eager-rw-penalty-bound}
    Let \(B\) be an invariant edge process on \(T_d\). Suppose its marginal distribution \(X\) has \(\EE(X)=0\) and \(\sigma^2=\Var(X)\lt\infty\). Then \[
        \EE_\nu(B_{o,x_{-1}\wedge \xi})\le C\frac{\sigma}{p}
    \]
    where \(C=\frac{1}{\sqrt{2}-1}\)
\end{proposition}
\begin{proof}
    For every vertex \(x\in S_n(o)\), write the edges of the path \([o,x]\) as \(e_1(x),...,e_n(x)\). Then \[
        B_{o,x_{-1}\wedge \xi}=\sum_{\ell=1}^nB(e_\ell(x_{-1}))\one_{|x_{-1}\wedge \xi|\ge \ell}
    \]
    Now by the Cauchy-Schwartz inequality, \begin{equation}
        \label{eqn:eager-rw-penalty-bound-1}
        \begin{aligned}
            \EE_\nu(B_{o,x_{-1}\wedge \xi})&=\sum_{\ell=1}^n \EE_\nu\big(B(e_\ell(x_{-1}))\one_{|x_{-1}\wedge \xi|\ge \ell}\big) \\
            &\le \sum_{\ell=1}^n \sqrt{\EE_\nu\big(B^2(e_\ell(x_{-1}))\big)}\sqrt{\PP_\nu(|x_{-1}\wedge \xi|\ge \ell)}
        \end{aligned}
    \end{equation}
    By \Cref{lemma:eager-rw-backtrack-exponential-decay}, the second term satisfies \begin{equation}
    \label{eqn:eager-rw-penalty-bound-2}
        \sqrt{\PP_\nu(|x_{-1}\wedge \xi|\ge \ell)}\le \frac{1}{\sqrt{p}}\frac{1}{(d-1)^{\ell/2}}
    \end{equation}
    For the first term, note that for each edge \(e\) with \(|e|=\ell-1\) (i.e. one of its endpoints is in \(S_{\ell-1}(o)\) and the other in \(S_{\ell}(o)\)), there are exactly \((d-1)^{n-\ell}=d_n/d_\ell\) vertices \(x\in S_n(o)\) for which \(e_\ell(x)=e\). Hence  
    \begin{align*}
        \nu_\omega(o)\EE(B^2(e_\ell(x_{-1}))\mid\omega)&=\sum_{x\in S_n(o)}c_\omega(o,x)B^2(e_\ell(x)) \\
        &\le \left(1+\frac{\eeps}{d_n}\right)\sum_{x\in S_n(o)}B^2(e_\ell(x)) \\
        &=\left(1+\frac{\eeps}{d_n}\right)\cdot \frac{d_n}{d_\ell}\sum_{|e|=\ell-1}B^2(e) =\frac{d_n+\eeps}{d_\ell}\sum_{|e|=\ell-1}B^2(e)
    \end{align*}
    Therefore \begin{equation}
    \label{eqn:eager-rw-penalty-bound-3}
        \EE_\nu(B^2(e_\ell(x_{-1})))\le \frac{d_n+\eeps}{\EE(\nu(o))}\EE\big(B^2(e)\big)\le \frac{1}{p}\sigma^2
    \end{equation}
    where \(e\) is an arbitrary fixed edge. Putting \eqref{eqn:eager-rw-penalty-bound-2} and \eqref{eqn:eager-rw-penalty-bound-3} into \eqref{eqn:eager-rw-penalty-bound-1}, we obtain \[
        \EE_\nu(B_{o,{x_{-1}\wedge \xi}})\le \sum_{\ell=1}^n \frac{\sigma}{\sqrt{p}}\frac{1}{\sqrt{p}}\frac{1}{(d-1)^{\ell/2}}=\frac{\sigma}{p}\sum_{\ell=1}^{n} \frac{1}{(d-1)^{\ell/2}}\le C\frac{\sigma}{p}
    \]
    where we chose \[
        C=\sum_{\ell=1}^\infty\frac{1}{(\sqrt{2})^\ell}=\frac{1}{\sqrt{2}-1} \qedhere
    \]
\end{proof}
\vspace{1em}

\subsection{Proof of \Cref{thm:exceptional-rays-exist}}
Let \(A\) be a non-constant, ergodic edge process, whose marginal distribution \(X\) has finite variance. Without loss of generality, \(\EE(X)=0\). \Cref{thm:exceptional-rays-exist} claims that a.s. \[
    \ssust(A)\ge \delta_X
\]
where \(\delta_X\gt 0\) depends only on the law of \(X\). \\

We start with a reduction to \(\Aut_+(T_d)\)-ergodic processes: let \(A\) be an \(\Aut(T_d)\)-ergodic process, with law \(\mu\) and marginal distribution \(X\). Write \(\mu=\frac{1}{2}\mu_1+\frac{1}{2}\mu_2\) as in \Cref{prop:equal-mixture}. Assuming the theorem holds for \(\Aut_+(T_d)\)-ergodic processes, we apply it to \(\mu_1,\mu_2\). Since \(\mu_i\) share the same marginal distribution \(X\), we have \[
    \mu_i( \ssust \ge \delta_X )=1
\]
where \(\delta_X\) is as in the theorem. As a result, \(\mu(\ssust\ge \delta_X)=1\). \\

Now, borrowing inspiration from the proof of \Cref{prop:exceptional-rays-haggstrom-argument}, let \(A_1=A,A_2,A_3...\) be i.i.d. copies of \(A\), and consider the pointwise average \[
    B=B_k=\frac{1}{k}(A_1+...+A_k)
\]
Let \(a\) be another parameter representing a ``threshold'', above which paths are considered preferable. Define \(C=C_{k,a}\) by \begin{equation}
    \label{eqn:bounded-case-eager-rw-condition}
    C_{k,a}(e)=\one_{B(e)\ge a}  
\end{equation}
and let \(\xx=\xx^{(n,k,a,\eeps)}=\xx^{(n,\eeps,C_{k,a})}\) be the corresponding eager random walk. The remainder of this subsection is devoted to proving the following theorem (which, by the above, implies \Cref{thm:exceptional-rays-exist}).

\begin{theorem} \label{thm:eager-rw-achieves-exceptional-ray}
    Let \(A\) be a non-constant, \(\Aut_+(T_d)\)-ergodic edge process. Suppose its marginal distribution \(X\) satisfies \(\EE(X)=0\) and \(\sigma^2=\Var(X)\lt\infty\).
    There exist parameters \(n,k,a,\eeps\) such that the eager random walk \(\xx=\xx^{(n,k,a,\eeps)}\) a.s. satisfies \[
        \MM(\xi(\xx);A)\ge \delta_X\quad
    \]
    where \(\delta_X\gt 0\) depends only on the law of \(X\).
\end{theorem}

First, the walk \(\xx\) is a.s. transient, and in particular end-convergent. In addition, by the Howe-Moore property (\Cref{thm:howe-moore-property}), \(A\) is mixing, hence \(B\) must also be mixing, and in particular ergodic.

We also assume that \(n\) is even: in this case, we have \(\Aut(G_n)=\Aut_+(T_d)\), hence the walk \(\xx\) must be stationary. Since \(B\) is ergodic, \(\xx\) must also be ergodic by \Cref{thm:rwrn-ergodicity-criterion} (For odd \(n\), the walk oscillates between the tree's bipartite classes \(P,Q\); since \(B\) is only guaranteed to be \(\Aut_+(T_d)\)-invariant, the distribution of the environment seen from the particle might differ between them). \\

Next, we connect between the asymptotic means of \(A\) and \(B\) along \(\xi=\xi(\xx)\).

\begin{observation} \label{observation:mean-of-A-initial-estimate}
    Suppose \(\speed(\xx)\gt 0\). Then a.s.,
    \begin{equation} \label{eqn:mean-of-A-initial-estimate}
        \MM(\xi;A)= \MM(\xi;B)=\frac{\EE_\nu(B_{o,x_1})-2\EE_\nu(B_{o,x_{-1}\wedge \xi})}{\speed(\xx)}   
    \end{equation}
\end{observation}
\begin{proof}
    By linearity of limits, \[
        \MM(\xi;B)=\frac{1}{k}\sum_{i=1}^k\MM(\xi;A_i)=\MM(\xi;A)
    \]
    where the second equality follows since \(A\) and the walk \(\xx\) are ergodic. The second equality follows from \Cref{thm:mean-along-limiting-end}.
\end{proof}

In particular, given that \(\speed(\xx)\gt 0\), \(\MM(\xi;A)\) is positive if and only if the numerator on the right side of \eqref{eqn:mean-of-A-initial-estimate} is positive. We proceed by estimating both of its terms.

\begin{proposition} \label{prop:density-of-single-step}
\[
    \EE_\nu(B_{o,x_1})\ge a\left(1-\frac{\eeps}{\EE(\nu(o))}\right)
\]
\end{proposition}
\begin{proof}
    We have \begin{align*}
        \nu_\omega(o)\EE(B_{o,x_1}\mid \omega)&=\sum_{x\in S_n(o)}c_\omega(o,x)B_{o,x}  \\
        &=\frac{\eeps}{d_n}\sum_{x\in S_n(o)}B_{o,x}+\sum_{x\in S_n(o)}B_{o,x}\one_{(B_k)_{o,x}\ge a} \\
        &\ge \frac{\eeps}{d_n}\sum_{x\in S_n(o)}B_{o,x}+a\deg_{C}(o)
    \end{align*}
    Therefore \begin{align*}
        \EE_\nu(B_{o,x_1})&=\frac{\EE(\nu_\omega(o)\EE(B_{o,x_1}\mid \omega))}{\EE(\nu(o))} \\
        &\ge \frac{\eeps}{\EE(\nu(o))}\EE(B_{o,x})+\frac{a\EE(\deg_{C}(o))}{\EE(\nu(o))}=a\left(1-\frac{\eeps}{\EE(\nu)}\right)   
    \end{align*}
    Where the last equality follows since \(\EE(B_{o,x})=0\) and \(\nu(o)=\deg_C(o)+\eeps\).
\end{proof}

For the second term, let \[
    p=\PP(B_{x,y}\ge a)
\]
be the marginal probability of \(C_{k,a}\). Note that the standard deviation of the marginal of \(B\) is \(\sigma/\sqrt{k}\). Hence, by \Cref{prop:eager-rw-penalty-bound}, \[
    \EE_\nu(B_{o,x_{-1}\wedge \xi})\le \frac{C\sigma}{p\sqrt{k}}
\]
Together with \Cref{prop:density-of-single-step}, we obtain  \[
    \EE_\nu(B_{o,x_1})-2\EE_\nu(B_{o,x_{-1}\wedge \xi})\ge a\left(1-\frac{\eeps}{\EE(\nu(o))}\right)-2\frac{C\sigma}{p\sqrt{k}}
\]
Since \(\eeps\gt 0\) is arbitrary, it suffices to find parameters \(n,k,a\) such that \[
    \Delta:=a-2\frac{C\sigma}{p\sqrt{k}}\gt 0
\]

Let \(\theta\gt 0\) be such that \(\PP(Z\ge \theta)=0.45\), where \(Z\sim N(0,1)\). Fix a pair of vertices \(x,y\) with \(d_T(x,y)=n\). Let \(S_n\overset{d}{=}A_{x,y}\) be distributed as the sum along \([x,y]\), and let \(\sigma_n\) be its standard deviation. Choose \begin{equation} \label{eqn:theorem-1.4-choice-of-a}
    a=a_{n,k}=\frac{\theta\sigma_n}{\sqrt{k}}
\end{equation}

Then our estimate becomes \[
    \Delta=\frac{\theta\sigma_n}{\sqrt{k}}-\frac{2C\sigma}{p\sqrt{k}}=\frac{\theta}{\sqrt{k}}\left(\sigma_n-\frac{2C}{p\theta }\sigma\right)
\]
Now as \(k\to\infty\), the central limit theorem implies that \[
    p=\PP(B_{x,y}\ge a)=\PP\left(\frac{(A_1)_{x,y}+...+(A_k)_{x,y}}{\sigma_n \sqrt{k}}\ge \theta\right)\xrightarrow[k\to\infty]{} \PP(Z\ge \theta)=0.45
\]
Hence for \(k\) large enough, \(p\ge 0.4\). This implies, by \Cref{corollary:eager-rw-speed-lower-bound}, that a.s. \[
    \speed(\xx)\ge n-\frac{2}{p}\ge n-5
\]
and for \(n\ge 6\) the walk is guaranteed to have positive speed (such that \Cref{observation:mean-of-A-initial-estimate}, and our estimate \(\Delta\) derived from it, apply). Therefore, it suffices to show we can choose \(n\) large enough, independently of \(k\), such that, \[
    \sigma_n\ge 5\frac{C}{\theta}\sigma
\]
The next result says a little more:

\begin{proposition} \label{prop:variance-growth}
    Let \(S_n\) be distributed as the sum of a path of length \(n\) with respect to \(A\). Then 
    \begin{itemize}
        \item For \(d\ge 4\), \(\Var(S_{2^n})\ge \left(2\left(\frac{d-2}{d-1}\right)\right)^n\Var(X)\ge \left(\frac{4}{3}\right)^n\Var(X)\)
        \item For \(d=3\), \(\Var(S_{2\cdot3^n})\ge \Var(S_{3^n})\ge \left(\frac{5}{4}\right)^n\Var(X)\)
    \end{itemize}
\end{proposition}

\begin{proof}
    For \(d\ge 4\), consider a path of length \(2n\) in the tree. Let \(S_n',S_n''\) be the sums along its first and second halves, respectively, and let \[
        \rho_0:=\rho(S_n',S_n'')=\frac{\Cov(S_n',S_n'')}{\Var(S_n)}
    \]
    be their correlation coefficient. Then \[
        \Var(S_{2n})=\Var(S_n')+\Var(S_n'')+2\Cov(S_n',S_n'')=(2+2\rho_0)\Var(S_n)
    \]
    To bound \(\rho_0\) from below, note that for a given vertex \(v\) in \(T_d\), we can pick \(d\) edge-disjoint paths of length \(n\) that start at \(v\). Let \(S_n^{(1)},...,S_n^{(d)}\) be the sums along these paths, then \[
        0\le \Var\left(S_n^{(1)}+...+S_n^{(d)}\right)=d\Var(S_n)+d(d-1)\rho_0\Var(S_n)
    \]
    which implies \(\rho_0\ge -\frac{1}{d-1}\). Therefore for each \(n\), \[
        \Var(S_{2n})\ge 2\left(1-\frac{1}{d-1}\right)\Var(S_n)
    \]
    and the result follows by induction on \(n\).

    For \(d=3\), the same argument gives that \(\rho_{0}\ge -\frac{1}{2}\), which gives only \(\Var(S_{2n})\ge \Var(S_n)\). To get our result we have to work a little harder: this time, consider a path of length \(3n\), split into three segments of length \(n\), with sums \(S',S'',S'''\). Then \begin{equation} \label{eqn:bounded-case-variance-growth-T3}
        \Var(S_{3n})=(3+4\rho_0+2\rho_1)\Var(S_n)
    \end{equation} 
    where \(\rho_0=\rho(S',S'')=\rho(S'',S''')\) is as before, and \(\rho_1:=\rho\left(S',S'''\right)\). Now, fix an \(n\)-path in \(T_3\). From each of its ends, we may pick two additional paths, such that all five paths are edge-disjoint. Denote that sums along those paths by \(S_n^{(1)},...,S_n^{(5)}\). Then \[
        0\le \Var\left(S_n^{(1)}+...+S_n^{(5)}\right)=(5+12\rho_0+8\rho_1)\Var(S_n)
    \] 
    Hence \(5+12\rho_0+8\rho_1\ge 0\), and simplifying we get \(3\rho_0+2\rho_1\ge -\frac{5}{4}\). But \(\rho_0\ge -\frac{1}{2}\) by the above, so \[
        3+4\rho_0+2\rho_1\ge 3-\frac{1}{2}-\frac{5}{4}\ge \frac{5}{4}
    \]
    
    Putting this back into \eqref{eqn:bounded-case-variance-growth-T3}, we conclude \[
        \Var(S_{3n})\ge \frac{5}{4}\Var(S_n)
    \]
    and again the result follows by induction.
\end{proof}

\begin{remark} \label{remark:variance-growth-restatement}
    \Cref{prop:variance-growth} may be restated as follows: there exists a global parameter \(\alpha=\alpha_d\gt 0\), depending only on \(d\), such that for every non-constant invariant edge process \(A\) on \(T_d\), \[
        \limsup_{n\to\infty} n^{-\alpha}\Var(S_n)\gt 0
    \]
    More precisely, for \(d\ge 4\) we take \(\alpha_d=\log_2\left(2\frac{d-2}{d-1}\right)\); then along the indices \(m_n=2^n\), we have \[
        m_n^{-\alpha_d}\Var(S_{m_n})\ge \left(2\left(\frac{d-1}{d-2}\right)\right)^{-n}\Var(S_{2^n})\ge \Var(X)\gt 0
    \]
    Similarly, for \(d=3\) we may take \(\alpha_3=\log_3\left(\frac{5}{4}\right)\) and \(m_n=3^n\).
\end{remark}

\begin{remark}
    Let \(\gamma\) be a ray in \(T_d\), and \((S_n)_n\) the partial sums of the first \(n\) edges along \(\gamma\), as above. \Cref{remark:variance-growth-restatement} implies, in particular, for any non-constant, invariant edge process, \[
    \limsup_{n}\Var(S_n)=\infty
    \]
    For stationary processes on \(\NN\) or \(\ZZ\), this is known to be false: the central example is the family of \textit{coboundary} processes, which are of the form \(X_n=Y_{n+1}-Y_n\) for some other stationary process \(Y=(Y_n)_n\) (with finite variance). Indeed, for such a process, \[
    S_n=\sum_{i=0}^{n-1}X_i=Y_n-Y_0
    \]
    which implies that \(\Var(S_n)=\Var(Y_n-Y_0)\le 4\Var(Y_0)\lt \infty\). In fact, it is well known that under mild conditions (namely, that \(\Cov(X_0,X_n)\underset{n\to\infty}{\longrightarrow} 0\)), \(\sup_n \Var(S_n)\lt \infty\) if and only if the process \(X\) is of such form (see \cite[Theorem 18.2.2]{ibragimov1971independent}). \\
    A corollary of \Cref{prop:variance-growth}, then, is that there are no non-constant invariant edge processes on \(T_d\) whose restriction to a ray is a coboundary process.
\end{remark}
\vspace{1em}
With \Cref{prop:variance-growth} established, the proof of \Cref{thm:eager-rw-achieves-exceptional-ray} is almost complete: choose \(n\) even and large enough such that \(\speed(\xx)\gt 0\) and \begin{equation}
    \sigma_n\ge \frac{5C}{\theta}\sigma
\end{equation}
(\(\theta\gt \frac{1}{8}\) and \(C\lt 3\), so it suffices that \(\sigma_n\ge 120\sigma\); by \Cref{prop:variance-growth}, this translates to \(n\ge2^{34}\) for \(d\ge 4\), or \(n\ge 3^{43}\) for \(d=3\)). Then for \(a=a_{n,k}\) as in \eqref{eqn:theorem-1.4-choice-of-a}, and for large enough \(k\), \begin{equation} \label{eqn:Delta-estimate}
    \Delta\ge \frac{\theta}{\sqrt{k}}\left(\sigma_n-\frac{5C}{\theta }\sigma\right)\gt  0
\end{equation}
For \(\eeps\) small enough, we obtain \[
    \MM(\xi;A)\ge \frac{1}{2}\frac{\Delta}{\speed(\xx)}\gt 0
\]
To finish the proof, we need to find some \(\delta_X\gt 0\), depending only on the marginal distribution \(X\), such that \(\MM(\xi;A)\ge \delta_X\). Note that \(n\) may be chosen globally, that is, independently of the process \(A\); \(a\) is determined by \(n\) and \(k\) and \(\eeps\) is arbitrary, so we are left with showing that \(k\) can be taken depending only on the law of \(X\). This is not immediate -- since we require \(k\) to be large enough that \(p=\PP(B_{x,y}\ge a)\approx \PP(Z\ge \theta)\), apriori we can only guarantee that \(k\) depends on the marginal distribution of an \(n\)-path, that is on \(S_n\). \\

If \(X\) has a finite (absolute) third moment \(\rho:=\EE(|X|^3)\lt\infty\), we may bound the difference \(|p-\PP(Z\ge \theta)|\) using the Berry-Esseen theorem (see, for example, \cite[Theorem 3.4.17]{durrett2019}):

\begin{theorem}[Berry-Esseen]
    Let \(Y_1,Y_2,...,Y_k\) be i.i.d. with \(\EE(Y_i)=0\), \(\EE(Y_i^2)=\sigma^2\lt \infty\) and \(\rho:=\EE(|Y_i|^3)\lt \infty\). Then for every \(\theta\), \[
        \left|\PP\left(\frac{Y_1+...+Y_k}{\sigma\sqrt{k}}\le \theta\right)-\PP(Z\le \theta)\right|\le \frac{C\rho}{\sigma^3\sqrt{k}}
    \]
    Where \(Z\sim N(0,1)\) and \(C\) is a global constant.
\end{theorem}

We apply the theorem to \(S_n\): after picking \(n\) appropriately (and globally), we have \(\sigma_n\ge n^{\alpha/2}\sigma\), so the error term in Berry-Esseen's theorem is bounded by \[
    \frac{C\EE(|S_n|^3)}{\sigma_n^3\sqrt{k}}\le \frac{Cn^3\EE(|X|^3)}{n^{3\alpha/2}\sigma\sqrt{k}}
\]
So \(p\ge 0.4\) is guaranteed whenever \[
    \frac{Cn^3\EE(|X|^3)}{n^{3\alpha/2}\sigma\sqrt{k}}\le 0.05\Longleftrightarrow 
    k\ge \frac{400C^2n^{6-3\alpha}\EE(|X|^3)}{\Var(X)}
\]

If \(X\) does not have a finite third moment, Berry-Esseen does not apply; we use instead the following, finer result due to Feller (see \cite{Feller1968}).

\begin{theorem}[Feller] \label{thm:feller}
    Let \(Y_1,Y_2,...\) be i.i.d. with \(\EE(Y_i)=0\), \mbox{\(\EE(Y_i^2)=\sigma^2\lt\infty\)}. Then for every \(M\) and \(\theta\), \[
        \left|\PP\left(\frac{Y_1+...+Y_k}{\sigma\sqrt{k}}\le \theta\right)-\PP(Z\le \theta)\right|\le C\left(\frac{\EE(Y_i^2\one_{|Y_i|\gt M})}{\sigma^2}+\frac{\EE(|Y_i|^3\one_{|Y_i|\le M})}{\sigma^3\sqrt{k}}\right)
    \]
    Where \(Z\sim N(0,1)\) and \(C\) is a global constant.
\end{theorem}

Again, we apply the theorem to \(S_n\). We wish to show that the bound can be made arbitrarily small using only the distribution of \(X\). The second term is bounded by \[
    \frac{\EE(|S_n|^3\one_{|S_n|\le M})}{\sigma_n^3\sqrt{k}}\le \frac{M\sigma_n^2}{\sigma_n^3\sqrt{k}}\le\frac{M}{n^{\alpha/2}\sigma\sqrt{k}}
\]
For the first term, write \(S_n=X_1+...+X_n\) as the sum of \(n\) random variables \(X_i\overset{d}{=} X\). To bound the expected value \(\EE(Y_i^2\one_{|Y_i|\gt M})\) using the distribution of a single element, we introduce another truncation, this time on the summands \(X_i\): for a fixed \(L\gt 0\), write \(X_i=X_i'+X_i''\) where \[
    X_i'=X_i\one_{|X_i|\le L},\quad X_i''=X_i\one_{|X_i|\gt L},
\] 
and \(S_n'=\sum X_i',S_n''=\sum X_i''\). Recall that for any real numbers \(a_1,...,a_m\), \begin{equation}
    \label{eqn:square-of-sum-bound}
    \left(\sum_{i=1}^m a_i\right)^2\le m\sum_{i=1}^m a_i^2  
\end{equation}
(this follows from either Cauchy-Schwartz or Jensen's inequality). In particular for \(m=2\), \[
    \EE(S_n^2\one_{|S_n|\gt M})\le 2\EE((S_n')^2\one_{|S_n|\gt M})+2\EE((S_n'')^2\one_{|S_n|\gt M})
\]
For the first term, by Chebychev's inequality, \[
    \EE((S_n')^2\one_{|S_n|\gt M})\le (Ln)^2\PP(|S_n|\gt M)\le \frac{L^2n^2\Var(S_n)}{M^2}\le \frac{L^2n^4\sigma^2}{M^2}
\]
For the second, using \eqref{eqn:square-of-sum-bound} again, \begin{align*}
    \EE((S_n'')^2\one_{|S_n|\gt M})
    \le n\EE\left(\sum_{i=1}^n (X_i'')^2\one_{|S_n|\gt M}\right)
    &\le n^2\EE((X_i'')^2) \\
    &=n^2\EE(X^2\one_{|X|\gt L})
\end{align*}

Now, let \(\eta=0.05/C\). Choose the parameters \(L,M,k\) as such. \begin{itemize}
    \item Choose \(L\) such that \(n^2\EE(X^2\one_{|X|\gt L})\lt \frac{\eta}{3}\).
    \item Then, choose \(M\) such that \(\frac{L^2n^4\sigma^2}{M^2}\le \frac{\eta}{3}\).
    \item Then, choose \(k\) such that \(\frac{M}{n^{\alpha/2}\sigma\sqrt{k}}\le \frac{\eta}{3}\).
\end{itemize}

By Feller's theorem and the computations above, we obtain that for these choices, \[
    |p-\PP(Z\ge \theta)|\le C\eta= 0.05
\]
and in particular \(p\ge 0.4\), as desired. Finally, note that the choice of \(L,M,k\) depends only on the marginal distribution \(X\) (and the global parameter \(n\)); plugging this choice of \(k\) into \eqref{eqn:Delta-estimate}, together with our choice of \(n\) and a small enough \(\eeps\), yields a gap \[
    \MM(\xi;A)\ge \delta_X\gt 0
\]
that depends only on \(X\), and the proof is complete.  \\

\subsection{Proof of \Cref{thm:vertex-exceptional-rays-exist}}
We now prove the analog of \Cref{thm:exceptional-rays-exist} for vertex processes. Recall that \(T_d\) is bipartite with partition \((P,Q)\) (without loss of generality, \(o\in P\)). For \(a,b\in \RR\), denote by \(\omega_{a,b}\) the vertex configuration on \(T_d\) defined by \(\omega_{a,b}|_P\equiv a\) and \(\omega_{a,b}|_Q\equiv b\). As in the introduction, consider the family of vertex processes \((\mathcal{S}_{a,b})_{a,b}\) distributed as \[
    \operatorname{Law}(\mathcal{S}_{a,b})=\tfrac{1}{2}\delta_{\omega_{a,b}}+\tfrac{1}{2}\delta_{\omega_{b,a}}
\]
\Cref{thm:vertex-exceptional-rays-exist} claims that among ergodic processes with finite marginal variance, these are the only ones satisfying \(\ssust(A)=\EE(X)\). We obtain it by an appropriate reduction to edge processes: let \(A\) be an ergodic vertex process on \(T_d\), whose marginal distribution \(X\) has finite variance. Define an edge process \(\tilde{A}\) by \[
    \tilde{A}(x,y):=\tfrac{1}{2}(A(x)+A(y))
\]
Note that \(\tilde{A}\) is ergodic, and its marginal \(\tilde{X}\) satisfies \(\EE(\tilde{X})=\EE(X)\) and \(\Var(\tilde{X})\le \Var(X)\lt\infty\). \Cref{thm:vertex-exceptional-rays-exist} follows from the combination of the following two propositions.
\begin{enumerate}
    \item If \(\tilde{X}\) is a.s. constant, then \( A\overset{d}{=}\mathcal{S}_{a,b}\) for some \(a,b\).
    \item Otherwise, \(\ssust(A)\gt \EE(X)\).
\end{enumerate}

\begin{proposition} \label{prop:vertex-processes-first-part}
    Suppose that \(\tilde X\) is a.s. constant. Then \(A\overset{d}{=}\mathcal{S}_{a,b}\) for some \(a,b\).
\end{proposition}

\begin{proof}
    Let \(\mu\) be the law of \(A\). For every vertex configuration \(\omega\), let \[
        f(\omega)(x,y)=\tfrac{1}{2}(\omega(x)+\omega(y))
    \]
    be the corresponding edge configuration, such that the law of \(\tilde A\) is \(\tilde \mu:=f_{\ast}\mu\). Our assumption that \(\tilde{X}\) is a.s. constant means that \(\tilde \mu\) is supported on constant edge configurations \(\tilde\omega\equiv c\). Given such a configuration \(f(\omega)\equiv c\), let \(a:=\omega(o)\) and \(b:=2c-a\). By assumption for every edge \((x,y)\), \[
        c=f(\omega)(x,y)=\tfrac{1}{2}(\omega(x)+\omega(y))
    \]
    hence \(\omega(x)=a\) if and only if \(\omega(y)=b\). In particular the neighbors \(x\sim o\) must all have \(\omega(x)=b\), their neighbors must all have a value of \(a\), and so on -- by induction, we get that \(\omega|_P\equiv a\) and \(\omega|_Q\equiv b\), that is \(\omega=\omega_{a,b}\). In other words, because \(\tilde\mu\) is supported on constant edge configurations, \(\mu\) must be supported on \(\Omega'=\set{\omega_{a,b}\mid a,b\in \RR}\). \\ 

    Now, consider the random variables \begin{align*}
        \alpha&=\alpha(\omega):=\inf_{v\in V(T_d)}\omega(v) \\       
        \beta&=\beta(\omega):=\sup_{v\in V(T_d)}\omega(v)       
    \end{align*}
    These are automorphism-invariant and \(\mu\) is ergodic, hence \(\alpha\equiv a, \beta\equiv b\) are \(\mu\)-a.s. constant. On the other hand \(\mu\) is supported on \(\Omega'\); taking the intersection we conclude that \(\mu\) is supported on \[
        \set{\omega\in \Omega'\mid \alpha(\omega)=a,\beta(\omega)=b}=\set{\omega_{a,b},\omega_{b,a}}
    \]
    or in other words, \[
        \mu=p_1\delta_{\omega_{a,b}}+p_2\delta_{\omega_{b,a}}
    \]
    Finally since \(A\) is \(\Aut(T_d)\)-invariant, and there exists an automorphism mapping the two bipartite classes \(P,Q\) onto one another, we must have \mbox{\(p_1=p_2=\frac{1}{2}\)}. Thus \[
        \mu=\tfrac{1}{2}\delta_{\omega_{a,b}}+\tfrac{1}{2}\delta_{\omega_{b,a}}=\operatorname{Law}(\mathcal{S}_{a,b}) \qedhere
    \]
\end{proof} 
\vspace{1em}

For the second point, by \Cref{thm:exceptional-rays-exist}, \(\ssust(\tilde A)\gt \EE(\tilde X)=\EE(X)\). If the process \(A\) is bounded, i.e. \(|X|\le M\) a.s. for some \(M\), then for any finite path \((\gamma_0,...,\gamma_n)\), \begin{equation} \label{eqn:bounded-vertex-processes-uniform-bound}
    \frac{1}{n}|A_{o,\gamma_n}-\tilde A_{o,\gamma_n}|=\frac{1}{2n}|A(o)+A(\gamma_n)|\le \frac{M}{n}\xrightarrow[n\to\infty]{} 0 
\end{equation}
Thus a.s. for all ends \(\zeta\), \[
    \MM_\ast(\zeta;A)=\MM_\ast(\zeta;\tilde A)
\]
hence \(\ssust(A)=\ssust(\tilde A)\gt \EE(X)\).

For a general, unbounded process, we cannot hope to have a uniform bound as in \eqref{eqn:bounded-vertex-processes-uniform-bound}. However we can still guarantee that \(\MM_\ast(\zeta;A)=\MM_\ast(\zeta;\tilde A)\) for ``most'' ends. Call an end \(\zeta\) \textbf{tame} if for some (every) representative ray \(\gamma=[x,\zeta)\), \[
    \frac{|A(\gamma_n)|}{n}\xrightarrow[n\to\infty]{} 0
\]
Every tame end \(\zeta\) satisfies \(\MM_\ast(\zeta;A)=\MM_\ast(\zeta;\tilde A)\) (by an argument similar to \eqref{eqn:bounded-vertex-processes-uniform-bound}). In fact, non-tame ends form an exceptional set if they exist, since for a fixed end \(\zeta=[\gamma]\), by the ergodic theorem, \[
    \frac{1}{n}|A(\gamma_n)|=\frac{1}{n}|A|_{o,\gamma_n}-\frac{1}{n}|A|_{o,\gamma_{n-1}}\xrightarrow[n\to\infty]{\as} \EE(|X|)-\EE(|X|)=0
\]
Nevertheless, non-tame rays can exist. We provide an example for completeness.

\begin{example} \label{example:untame-rays}
    Let \(X\) be a random variable distributed as \[
        \PP(X=n^2)=2^{-n}
    \]
    Then \(\EE(X)=\sum n^22^{-n}=6\). Consider the i.i.d. vertex process \(A\) on \(T_3\) with marginal distribution \(X\). In the notation of \Cref{section:scloud-equals-sburst}, fix some directed edge \((x;y)\), and consider the set of vertices in \(T(x;y)\) at distance \(n\) from the root \(y\). For each such vertex \(z\), let \(\hat z\) denote its parent (i.e. the predecessor of \(z\) on the path \([y,z]\)). Let \(C_n(x;y)\) denote the set of edges \((\hat z,z)\) that have value at least \(n^2\). Then \[
        \PP(C_n(x;y)=\emptyset)=\PP\left(X\lt n^2\right)^{2^{n-1}}=\left(1-2^{-(n-1)}\right)^{2^{n}}\xrightarrow[n\to\infty]{} e^{-2}, 
    \]
    hence \(\PP(C_n(x;y)\neq \emptyset)=1-e^{-2}\). Thus, we may construct a (random) non-tame ray end using the same technique as in the proof of \Cref{corr:sburst-equals-scloud}: start with an arbitrary directed edge \(\vec{e}_0\). Given \(\vec{e}_n=(x_n;y_n)\), If \(C_n(\vec{e}_n)\neq \emptyset\), pick \(\vec{e}_{n+1}\) to be one of its edges uniformly; otherwise pick \(\vec{e}_{n+1}\) uniformly among all edges \((\hat{z},z)\) with \(z\in T(x_n;y_n)\cap S_n(y_n)\). Let \(\gamma\) be the ray obtained by this process.
    Since the process is i.i.d., a.s. an infinite number of the edges \(\vec{e}_n=(x_n;y_n)\) satisfy \(A(\vec{e}_n)\ge n^2\). On the other hand since the \(n\)-th step is of size \(n\), we have that \[
        d(y_n,y_0)=\sum_{j=1}^nd(y_j,y_{j-1})=\sum_{j=1}^n j=\frac{n(n+1)}{2}
    \]
    Hence the ray \(\gamma\) satisfies \[
        \limsup_{n\to\infty}\frac{A(\gamma_{n-1},\gamma_n)}{n}\ge \limsup_{n\to\infty} \frac{A(x_n;y_n)}{d(y_n,y_0)}\ge \lim_{n\to\infty}\frac{n^2}{\frac{1}{2}n(n+1)}=2
    \]
    and \(\gamma\) is a.s. not tame.
\end{example}

\begin{proposition} \label{prop:vertex-processes-second-part}
    Suppose that \(\tilde{X}\) is not a.s. constant. Then \(\ssust(A)\gt \EE(X)\)
\end{proposition}
\begin{proof}
    By \Cref{thm:eager-rw-achieves-exceptional-ray}, applied to the edge process \(\tilde{A}\), there exist parameters such that the eager random walk \(\xx\) a.s. satisfies \[
        \MM(\xi(\xx);\tilde{A})\ge \EE(\tilde{X})+\delta_{\tilde{X}}\gt \EE(X)
    \]
    We claim that the limiting end \(\xi=\xi(\xx)\) is a.s. tame with respect to \(A\) -- this means that \(\MM(\xi;A)=\MM(\xi;\tilde A)\gt \EE(X)\), as claimed. Define another edge process \(A'\) by \[
        A'(x,y)=\frac{1}{2}(|A(x)|+|A(y)|)
    \]
    In other words, \(A'\) is the result of applying the operation \(A\mapsto \tilde {A}\) to \(|A|\). In particular, \(A'\) is ergodic and has finite marginal variance. 
    Let \(\gamma=[o,\xi)\) be the walk's escape ray. By \Cref{thm:mean-along-limiting-end}, the limit \(\MM(\xi;A')\) a.s. exists. Therefore, \begin{align*}
        \frac{|A(\gamma_n)|}{n}&\le \frac{2A'(\gamma_{n-1},\gamma_n)}{n}=\frac{2}{n}A'_{o,\gamma_n}-\frac{2}{n}A'_{o,\gamma_{n-1}} \\
        &\xrightarrow[n\to\infty]{\as} 2\MM(\xi;A')-2\MM(\xi;A') = 0 
    \end{align*}
    as claimed.
\end{proof}

Together, \Cref{prop:vertex-processes-first-part} and \Cref{prop:vertex-processes-second-part} imply \Cref{thm:vertex-exceptional-rays-exist}.

Note that unlike the edge case, here even when \(A\) is not of the form \(\mathcal{S}_{a,b}\), we cannot find a gap \(\delta_X\), depending only on the distribution \(X\) of a single vertex, for which \(\ssust(A)\ge \EE(X)+\delta_X\):

\begin{example}
    Define a vertex percolation process \(A\) as such. Let \(S=\mathcal{S}_{0,1}\), and let \(B\) be i.i.d. vertex percolation on \(T_d\) with marginal probability \(\frac{1}{2}\). For each edge \(e\) independently, toss a coin with some small probability \(\eeps\) of success; if it succeeds take \(A(v):=B(v)\), and otherwise take \(A(v):=S(v)\).

    This process has a marginal distribution of \(\operatorname{Bernoulli}(1/2)\) for every \(\eeps\ge 0\). However, \[
        \scloud(A)\le \scloud(S+B)\le \scloud(S)+\scloud(B)=\frac{1}{2}+\scloud(B)
    \]
    and \(\scloud(B)\to 0\) as \(\eeps\to 0\) (for example, by \cite[Theorem 13.32]{lyonsperes}), hence the difference \(\ssust(A)-\EE(X)\) can be made arbitrarily small.
\end{example}

That being said, the translation to edge processes naturally yields a gap \(\delta_{\tilde X}\) that depends only on the distribution of \(\tilde X\), i.e. on the joint distribution of two adjacent vertices.
\endgroup

\paragraph*{Acknowledgments.} We thank Zemer Kosloff for useful comments. 
This research is supported by the Israeli Science Foundation grant no. 3191/25.

\paragraph*{Statement on AI Use.} During the preparation of this manuscript, Google Gemini was used to find an error in an earlier version of \Cref{thm:mean-along-limiting-end}, and to suggest the result by Feller (\Cref{thm:feller}) and the truncation argument that follows. All other ideas, results, and proofs are our own.

\printbibliography

\end{document}